\documentclass[aos,MSNbibl,seceqn,rotating,dvips]{arximspdf}
\usepackage{dcolumn}

\doi{10.1214/13-AOS1172} 
\volume{41}
\issue{6}
\pubyear{2013}
\firstpage{2877}
\lastpage{2904}

\makeatletter
\newcommand{\widebar}{\overline}
\newcolumntype{d}[1]{D{.}{.}{#1}}
\newtheorem{teo}{Theorem}
\newtheorem{lem}{Lemma}
\newtheorem{cor}{Corollary}
\newproclaim{exa}{Example}
\makeatother

\begin{document}
\begin{frontmatter}

\title{A general theory of particle filters in hidden Markov models
and some applications}
\runtitle{Particle filter}
\pdftitle{A general theory of particle filters in hidden Markov models and some applications}

\begin{aug}
\author[A]{\fnms{Hock Peng} \snm{Chan}\corref{}\ead[label=e1]{stachp@nus.edu.sg}\thanksref{t1}}
\and
\author[B]{\fnms{Tze Leung} \snm{Lai}\ead[label=e2]{lait@stat.stanford.edu}\thanksref{t2}}
\runauthor{H. P. Chan and T. L. Lai}
\affiliation{National University of Singapore and Stanford University}
\address[A]{Department of Statistics\\
\quad and Applied Probability\\
National University of Singapore\\
6 Science Drive 2\\
Singapore 117546\\
Singapore\\
\printead{e1}} 
\address[B]{Department of Statistics\\
Stanford University\\
Sequoia Hall\\
Stanford, California 94305-4065\\
USA\\
\printead{e2}}
\end{aug}
\thankstext{t1}{Suported by the National University of Singapore Grant
R-155-000-120-112.}
\thankstext{t2}{Supported by NSF Grant DMS-11-06535.}

\received{\smonth{5} \syear{2012}}
\revised{\smonth{9} \syear{2013}}

%
\begin{abstract}
By making use of martingale representations, we derive the asymptotic
normality of particle filters in hidden
Markov models and a relatively simple formula for their asymptotic
variances. Although repeated resamplings result
in complicated dependence among the sample paths, the asymptotic
variance formula and martingale representations
lead to consistent estimates of the standard errors of the particle
filter estimates of the hidden states.
\end{abstract}

%
\pdfkeywords{60F10, 65C05, 60J22, 60K35, Particle filter, importance sampling and resampling, standard error}

\end{frontmatter}

\section{Introduction}\label{sec1}

Let $\mathbf{X}= \{ X_t, t \geq1 \}$ be a Markov chain and let $Y_1,
Y_2, \ldots$ be conditionally independent given $\mathbf{X}$, such that
%
\begin{equation}
\label{hmm} X_t \sim p_t(\cdot| X_{t-1}),\qquad Y_t \sim g_t(\cdot|X_t),
\end{equation}
in which $p_t$ and $g_t$ are density functions with respect to some
measures $\nu_X$ and $\nu_Y$, and $p_1(\cdot|X_0)$ denotes the initial
density $p_1(\cdot)$ of $X_1$. Gordon, Salmond and Smith~\cite{GSS93}
introduced \textit{particle filters} for Monte Carlo estimation of
$E[\psi(X_T)|Y_1, \ldots, Y_T]$. More generally, letting
$\mathbf{X}_t=(X_1,\ldots,X_t)$, $\mathbf{Y}_t=(Y_1,\ldots,Y_t)$, and
$\psi$ be a measurable real-valued function of $\mathbf{X}_T$, we
consider estimation of $\psi_T:=E[\psi(\mathbf{X}_T)|\mathbf{Y}_T]$.
The density function of $\mathbf{X}_T$ conditional on $\mathbf{Y}_T$ in
the above hidden Markov model (HMM) is
%
\begin{equation}
\label{pt} \tilde p_T(\mathbf{x}_T|
\mathbf{Y}_T) \propto\prod_{t=1}^T
\bigl[p_t (x_t|x_{t-1} ) g_t(Y_t|x_t)
\bigr].
\end{equation}
However, this conditional distribution is often difficult to sample
from and the normalizing constant is also difficult to compute for
high-dimensional or complicated state spaces. Particle filters that use
sequential Monte Carlo (SMC) methods involving importance sampling and
resampling have been developed to circumvent this difficulty.
Asymptotic normality of the particle filter estimate of $\psi_T$ when
the number of simulated trajectories becomes infinite has also been
established \cite {Cho04,Del04,DM08,Kun05}. However, no explicit
formulas are available to estimate the standard error of $\hat\psi_T$
consistently.

In this paper, we provide a comprehensive theory of the SMC estimate
$\hat\psi_T$, which includes asymptotic normality and consistent
standard error estimation. The main results are stated in Theorems
\ref{thm3a}~and~\ref{thm4a} in Section~\ref{sec2} for the case in which
bootstrap resampling is used, and extensions are given in
Section~\ref{sec4} to residual Bernoulli (instead of bootstrap)
resampling. The proof of Theorem \ref{thm3a} for the case where
bootstrap resampling is used at every stage, given in
Section~\ref{sec3.2}, proceeds in two steps. First we assume that the
normalizing constants [i.e., constants of proportionality in
(\ref{pt})] are easily computable. We call this the ``basic
prototype,'' for which we derive in Section~\ref{sec3.1} martingale
representations of $m(\tilde\psi_T-\psi_T)$ for novel SMC estimates
$\tilde\psi_T$ that involve likelihood ratios and general resampling
weights. We first encountered this prototype in connection with rare
event simulation using SMC methods in \cite{CL11}. Although traditional
particle filters also use similar sequential importance sampling
procedures and resampling weights that are proportional to the
likelihood ratio statistics, the estimates of $\psi_T$ are coarser
averages that do not have this martingale property.
Section~\ref{sec3.2} gives an asymptotic analysis of $m(\hat\psi_T -
\psi_T)$ as $m \rightarrow\infty$ in this case, by making use of the
results in Section~\ref{sec3.1}.

In contrast to Gilks and Berzuini \cite{GB01} who use particle set
sizes $m_t$ which increase with $t$ and approach $\infty$ in a certain
way to guarantee consistency of their standard error estimate which
differs from ours, we use the same particle set size $m$ for every $t$,
where $m$ is the number of SMC trajectories. As noted in
page~132 of~\cite{GB01}, ``the mode of convergence of the theorem is not directly
relevant to the practical context'' in which $m_1 = m_2 = \cdots=m$. A
major reason why we are able to overcome this difficulty in developing
a consistent standard error estimate lies in the martingale
approximation for our SMC estimate $\hat\psi_T$. We note in this
connection that for the basic prototype, although both martingale
representations of $m(\tilde\psi_T - \psi_T)$ developed in
Section~\ref{sec3.1} can be used to prove the asymptotic normality of
$\tilde\psi_T$ as the number $m$ of SMC trajectories approaches
$\infty$, only one of them leads to an estimable expression for the
asymptotic variance. Further discussion of the main differences between
the traditional and our particle filters is given in
Section~\ref{sec5}.

\section{Standard error estimation for bootstrap filters}\label{sec2}

Since $Y_1, Y_2, \ldots$ are the observed data, we will treat them as
constants in the sequel. Let $q_t(\cdot|\mathbf{x}_{t-1})$ be a
conditional density function with respect to $\nu_X$ such that
$q_t(x_t|\mathbf{x}_{t-1}) > 0$ whenever $p_t(x_t|x_{t-1}) > 0$. In the
case $t=1$, $q_t(\cdot|\mathbf{x}_{t-1})=q_1(\cdot)$.

\subsection{Bootstrap resampling at every stage}\label{sec2.1}
Let
%
\begin{equation}
\label{wt} w_t(\mathbf{x}_t) = p_t(x_t|x_{t-1})
g_t(Y_t|x_t)/q_t(x_t|
\mathbf{x}_{t-1}).
\end{equation}
To estimate $\psi_T$, a particle filter first generates $m$
conditionally independent random variables $\widetilde{X}{}_t^{1},
\ldots, \widetilde{X}{}_t^{m}$ at stage $t$, with
$\widetilde{X}{}_t^{i}$ having density
$q_t(\cdot|\mathbf{X}_{t-1}^{i})$, to form
$\widetilde{\mathbf{X}}{}_t^{i}=(\mathbf{X}_{t-1}^{i},
\widetilde{X}{}_t^{i})$ and then use normalized resampling weights
%
\begin{equation}
\label{weights} W_t^i = w_t \bigl(\widetilde{\mathbf{X}}{}_t^i \bigr) \Big/ \sum
_{j=1}^m w_t \bigl(\widetilde{\mathbf{X}}{}_t^j \bigr)
\end{equation}
to\vspace*{1pt} draw $m$ sample paths $\mathbf{X}_t^{j}$, $1 \leq j
\leq m$, from $\{ \widetilde{\mathbf{X}}{}_t^{i}, 1 \leq i \leq m \}$.
Note that $W_t^i$ are~the importance weights attached to
$\widetilde{\mathbf{X}}{}_t^i$ and that after resampling the
$\mathbf{X}_t^i$ have equal weights. The following recursive algorithm
can be used to implement particle filters with bootstrap resampling.
It\vspace*{-1pt} generates not only the
$\widetilde{\mathbf{X}}{}_T^{i}$ but also the ancestral origins
$A_{T-1}^i$ that are used to compute the standard error estimates. It
also computes $\widetilde{H}{}_t^i$ and $H_t^i$ recursively, where
%
\begin{eqnarray}
\widetilde{H}{}_t^i & = & \bar
w_1 \cdots\bar w_t h_t \bigl(\widetilde{\mathbf{X}}{}_t^i \bigr),
\qquad H_t^{i}=
\bar w_1 \cdots\bar w_t h_t \bigl(
\mathbf{X}_t^{i} \bigr),
\nonumber\\[-9pt]\label{24a} \\[-9pt]
\bar w_k & = & \sum_{j=1}^m
w_k \bigl(\widetilde{\mathbf{X}}{}_k^j \bigr) /
m, \qquad h_t(\mathbf{x}_t) = 1 \Big/ \prod
_{k=1}^t w_k(\mathbf{x}_k).\nonumber
\end{eqnarray}

\textit{Initialization}: let $A_0^i =i$ and $H_0^i=1$ for all $1 \leq i
\leq m$.

\textit{Importance sampling at stage $t=1,\ldots,T$}: generate
conditionally independent $\widetilde{X}{}_t^i$ from $q_t
(\cdot|\mathbf{X}_{t-1}^i)$, $1 \leq i \leq m$.

\textit{Bootstrap resampling at stage $t=1,\ldots,T-1$}: generate
i.i.d. random variables $B_t^1, \ldots, B_t^m$ such that
$P(B_t^1=j)=W_t^j$ and let
\[
\bigl(\mathbf{X}_t^j, A_t^j,
H_t^j \bigr) = \bigl(\widetilde{\mathbf{X}}{}_t^{B_t^j},
A_{t-1}^{B_t^j}, \widetilde{H}{}_t^{B_t^j} \bigr)\qquad\mbox{for } 1 \leq j \leq m.
\]

The SMC estimate is defined by
%
\begin{equation}
\label{23a} \hat\psi_T = (m \bar w_T)^{-1}
\sum_{i=1}^m \psi \bigl(\widetilde{\mathbf{X}}{}_T^{i} \bigr) w_T \bigl(\widetilde{\mathbf{X}}{}_T^{i} \bigr).
\end{equation}
Let $E_{\tilde p_T}$ denote expectation with respect to the probability
measure under which $\mathbf{X}_T$ has density (\ref{pt}), and let
$E_q$ denote that under which $X_t|\mathbf{X}_{t-1}$ has the
conditional density function $q_t$ for all $1 \leq t \leq T$. In
Section~\ref{sec3.2} we prove the following theorem on the asymptotic
normality of $\hat\psi_T$ and the consistency of its standard error
estimate. Define
%
\begin{eqnarray}\label{etah}
\eta_t &=& E_q \Biggl[ \prod
_{k=1}^t w_k(\mathbf{X}_k)
\Biggr], \nonumber
\\
h_t^*(\mathbf{x}_t) &=& \eta_t \Big/ \prod_{k=1}^t w_k(\mathbf{x}_k),\qquad \zeta_t=\eta_{t-1}/\eta_t,
\\
\Gamma_t(\mathbf{x}_t) &=& \prod
_{k=1}^t \bigl[w_k(\mathbf
{x}_k)+w_k^2(\mathbf{x}_k)
\bigr],\nonumber
\end{eqnarray}
with the convention $\eta_0=1$ and $h_0^* \equiv\Gamma_0 \equiv1$. Let
%
\begin{equation}
\label{kappat} L_T(\mathbf{x}_T) = \tilde
p_T(\mathbf{x}_T|\mathbf{Y}_T) \Big/ \prod
_{t=1}^T q_t(x_t|
\mathbf{x}_{t-1}).
\end{equation}
Let $f_0 = 0$ and define for $1 \leq t \leq T$,
%
\begin{equation}
\label{25} f_t(\mathbf{x}_t) \bigl[=f_{t,T}(
\mathbf{x}_t) \bigr] = E_q \bigl\{ \bigl[\psi(\mathbf
{X}_T)-\psi_T \bigr] L_T(
\mathbf{X}_T)| \mathbf{X}_t=\mathbf{x}_t \bigr\}.
\end{equation}

\begin{teo} \label{thm3a} Assume that $\sigma^2 < \infty$, where
$\sigma^2 = \sum_{k=1}^{2T-1} \sigma_k^2$ and
%
\begin{eqnarray}\label{tsigk}
\sigma_{2t-1}^2 &=& E_q \bigl\{
\bigl[f_t^2(\mathbf{X}_t)-f_{t-1}^2(
\mathbf{X}_{t-1}) \bigr] h_{t-1}^*(\mathbf{X}_{t-1})
\bigr\},
\nonumber\\[-8pt]\\[-8pt]
\sigma_{2t}^2 &=& E_q
\bigl[f_t^2(\mathbf{X}_t) h_t^*(
\mathbf{X}_t) \bigr].\nonumber
\end{eqnarray}
Then $\sqrt{m}(\hat\psi_T-\psi_T) \Rightarrow N(0,\sigma^2)$. Moreover,
if $E_q \Gamma_T(\mathbf{X}_T) < \infty$, then
$\hat\sigma{}^2(\hat\psi_T) \stackrel{p}{\rightarrow} \sigma^2$, where
for any real number $\mu$,
%
\begin{equation}
\label{420} \hat\sigma{}^2(\mu) = m^{-1} \sum
_{j=1}^m \biggl( \sum_{i\dvtx A_{T-1}^i =j}
\frac{w_T(\widetilde{\mathbf{X}}{}_T^i)}{\bar w_T} \bigl[\psi \bigl(\widetilde{\mathbf{X}}{}_T^i
\bigr)-\mu \bigr] \biggr)^2.
\end{equation}
\end{teo}

\subsection{Extension to occasional resampling}\label{sec2.2}

Due to the computational cost of resampling and the variability
introduced by the resampling step, \textit{occasional resampling} has
been recommended, for example, by Liu (Chapter~3.4.4 of \cite{Liu08}), who
propose to resample at stage $k$ only when cv$_k$, the coefficient of
variation of the resampling weights at stage $k$, exceeds some
threshold. Let $\tau_1 < \tau_2 < \cdots< \tau_r$ be the successive
resampling times and let $\tilde\tau(k)$ be the most recent resampling
time before stage $k$, with $\tilde\tau(k)=0$ if no resampling has been
carried out before time $k$. Define
%
\begin{equation}
\label{421} v_k(\mathbf{x}_k) = \prod
_{t=\tilde\tau(k)+1}^k w_t(\mathbf{x}_t),
\qquad V_k^i = v_k \bigl(\widetilde{\mathbf{X}}{}_k^i \bigr) /(m \bar v_k),
\end{equation}
where $\bar v_k = m^{-1} \sum_{j=1}^m
v_k(\widetilde{\mathbf{X}}{}_k^j)$. Note that if bootstrap resampling
is carried out at stage $k$, then $V_k^i$, $1 \leq i \leq m$, are the
resampling weights. Define the SMC estimator
%
\begin{equation}
\label{422} \hat\psi_{\mathrm{OR}} = \sum_{i=1}^m
V_T^i \psi \bigl(\widetilde{\mathbf{X}}{}_T^i
\bigr).
\end{equation}
In Section~\ref{sec4.3}, we show that the resampling times, which
resample at stage~$k$ when cv$_k^2$ exceeds some threshold $c$, satisfy
%
\begin{equation}
\label{rstar} r \stackrel{p} {\rightarrow} r^*\quad\mbox{and}\quad \tau_s
\stackrel {p} {\rightarrow} \tau_s^*\qquad\mbox{for } 1 \leq s \leq r^*\mbox{ as } m \rightarrow \infty,
\end{equation}
for some nonrandom positive integers $r^*, \tau_1^*, \ldots$ and
$\tau_{r^*}^*$.

In Section~\ref{sec4.2} we prove the following extension of Theorem
\ref{thm3a} to $\hat\psi_{\mathrm{OR}}$ and give an explicit formula
(\ref{sigOR}), which is assumed to be finite, for the limiting
variance~$\sigma^2_{\mathrm{OR}}$ of
$\sqrt{m}(\hat\psi_{\mathrm{OR}}-\psi_T)$.

\begin{teo} \label{thm4a} Under (\ref{rstar}), $\sqrt {m}(\hat\psi_{\mathrm{OR}}-\psi_T) \Rightarrow N(0,\sigma^2_{\mathrm{OR}})$ as $m
\rightarrow\infty$. Moreover, if $E_q \Gamma_T (\mathbf{X}_T) <
\infty$, then $\hat\sigma{}^2_{\mathrm{OR}}(\hat\psi_{\mathrm{OR}})
\stackrel {p}{\rightarrow} \sigma^2_{\mathrm{OR}}$, where for any real
number $\mu$,
\[
\hat\sigma{}^2_{\mathrm{OR}}(\mu) =m^{-1} \sum
_{j=1}^m \biggl( \sum_{i\dvtx A_{\tau_r}^i =j}
\frac{v_T(\widetilde{\mathbf{X}}{}_T^i)}{\bar v_T} \bigl[\psi \bigl(\widetilde{\mathbf{X}}{}_T^i
\bigr)-\mu \bigr] \biggr)^2.
\]
\end{teo}

\subsection{A numerical study}\label{sec2.3}

Yao \cite{Yao84} has derived explicit formulas for $\psi_T =
E(X_T|\mathbf{Y}_T)$ in the normal mean shift model $\{ (X_t,Y_t)\dvtx
t \geq1 \}$ with the unobserved states $X_t$ generated recursively by
\[
X_t = \cases{X_{t-1}, &\quad with prob. $1-\rho$,
\cr
Z_t, &\quad with prob. $\rho$,}
\]
where $0 < \rho< 1$ is given and $Z_t \sim N(0,\xi)$ is independent of
$\mathbf{X}_{t-1}$. The observations are $Y_t=X_t+\varepsilon_t$, in
which the $\varepsilon_t$ are i.i.d. standard normal, $1 \leq t \leq
T$. Instead of $\mathbf{X}_T$, we simulate using the bootstrap filter
the change-point indicators $\mathbf{I}_T:= (I_1,\ldots, I_T)$, where
$I_t = \mathbf {1}_{\{ X_t \neq X_{t-1} \}}$ for $t \geq2$ and $I_1=1$;
this technique of simulating $E(X_T|\mathbf{Y}_T)$ via realizations of
$E(X_T|\mathbf {I}_T,\mathbf{Y}_T)$ is known as Rao--Blackwellization.
Let $C_t = \max\{ j \leq t\dvtx  I_j=1 \}$ denote the most recent
change-point up to time $t$, $\lambda_t=(t-C_t+1+\xi^{-1})^{-1}$ and
$\mu_t=\lambda_t \sum_{i=C_t}^t Y_i$. It can be shown that
$E(X_T|\mathbf{I}_T,\mathbf{Y}_T) = \mu _T$ and that for $t \geq2$, $P
\{ I_t=1 | \mathbf{I}_{t-1}, \mathbf{Y}_t \} = a(\rho )/\{
a(\rho)+b(\rho) \}$, where $a(\rho) = \rho\phi_{0,1+\xi}(Y_t)$ and
$b(\rho)= (1-\rho) \phi_{\mu_{t-1},1+\lambda_{t-1}} (Y_t)$, in which
$\phi_{\mu,V}$ denotes the $N(\mu,V)$ density function.

\begin{table}
\tabcolsep=0pt
 \caption{Simulated values of $\hat\psi_{\mathrm{OR}}
\pm\hat\sigma_{\mathrm{OR}} (\hat\psi_{\mathrm{OR}})/\sqrt{m}$. An
italicized value denotes that $\psi_T$ is outside the interval but
within two standard errors of $\hat\psi_{\mathrm{OR}}$; a~boldfaced
value denotes that $\psi_T$ is more than two standard errors away from
$\hat\psi_{\mathrm{OR}}$}\label{tab1}
{\fontsize{8.5pt}{11pt}\selectfont{\begin{tabular*}{\textwidth}{@{\extracolsep{\fill}}@{}lccccc@{}}
\hline $\bolds{c}$ & \multicolumn{1}{c}{$\bolds{T=200}$} &
\multicolumn{1}{c}{$\bolds{T=400}$}
            & \multicolumn{1}{c}{$\bolds{T=600}$} & \multicolumn{1}{c}{$\bolds{T=800}$}
            & \multicolumn{1}{c}{$\bolds{T=1000}$}
\\
\hline
$\infty$ & \textit{$-$0.3718${}\pm{}$0.0067} & \textit{$-$0.1791${}\pm{}$0.0281} & \textit{$-$0.3022${}\pm{}$0.0136} & 0.5708${}\pm{}$0.0380 & $-$0.5393${}\pm{}$0.0180 \\
50 & $-$0.3592${}\pm{}$0.0019 & $-$0.1453${}\pm{}$0.0056 &\textit{$-$0.2768${}\pm{}$0.0020} & 0.6004${}\pm{}$0.0030 & $-$0.5314${}\pm{}$0.0024 \\
10 & \textit{$-$0.3622${}\pm{}$0.0019} & $-$0.1421${}\pm{}$0.0032 & $-$0.2802${}\pm{}$0.0029 & 0.6013${}\pm{}$0.0024 & $-$0.5287${}\pm{}$0.0015 \\
\phantom{0}5 & \textit{$-$0.3558${}\pm{}$0.0018} &\textit{$-$0.1382${}\pm{}$0.0027} & $-$0.2777${}\pm{}$0.0026 & 0.6004${}\pm{}$0.0022 &\textit{$-$0.5321${}\pm{}$0.0012} \\
\phantom{0}2 & $-$0.3584${}\pm{}$0.0013 & $-$0.1421${}\pm{}$0.0024 & $-$0.2798${}\pm{}$0.0018 & 0.6030${}\pm{}$0.0020 & \textit{$-$0.5316${}\pm{}$0.0008} \\
\phantom{0}1& \textit{$-$0.3615${}\pm{}$0.0013} & $-$0.1402${}\pm{}$0.0027 &\textit{$-$0.2761${}\pm{}$0.0019} & \textit{0.5997${}\pm{}$0.0019} & $-$0.5291${}\pm{}$0.0014 \\
\phantom{0}0.5 & \textit{$-$0.3610${}\pm{}$0.0016} & $\bolds{-0.1488\pm 0.0036}$ &$-$0.2785${}\pm{}$0.0020 & \textit{0.6057${}\pm{}$0.0021} & $-$0.5315${}\pm{}$0.0014 \\
\phantom{0}0 &$-$0.3570${}\pm{}$0.0032 & \textit{$-$0.1485${}\pm{}$0.0067} & $-$0.2759${}\pm{}$0.0033 &0.6033${}\pm{}$0.0031 & $-$0.5301${}\pm{}$0.0047\\
\hline
$\psi_T$ & $-$0.3590 & $-$0.1415 & $-$0.2791 & 0.6021 & $-$0.5302 \\
\hline
\end{tabular*}}}
\end{table}

\begin{table}[b]
\tabcolsep=0pt
\tablewidth=305pt
\caption{Fraction of confidence intervals $\hat\psi_{\mathrm{OR}} \pm\hat\sigma/\sqrt{m}$ containing the true
mean $\psi_T$}\label{tab2}
\begin{tabular*}{305pt}{@{\extracolsep{\fill}}@{}lcccc@{}}
\hline
& \multicolumn{2}{c}{$\bolds{\hat\sigma{}^2=\hat\sigma{}^2_{\mathrm{OR}} (\hat\psi_{\mathrm{OR}})}$}
& \multicolumn{2}{c}{$\bolds{\hat\sigma{}^2 =\widehat V_T}$}
\\[-6pt]
& \multicolumn{2}{c}{\hrulefill} & \multicolumn{2}{c}{\hrulefill}
\\
$\bolds{T}$ & \textbf{1 se} & \textbf{2 se} & \textbf{1 se} & \textbf{2 se}\\
\hline
\phantom{0}200 & 0.644 & 0.956 & 0.980 & 1.000 \\
\phantom{0}400 & 0.652 & 0.948 & 0.998 & 1.000 \\
\phantom{0}600 & 0.674 & 0.958 & 1.000 & 1.000 \\
\phantom{0}800 & 0.716 & 0.974 & 1.000 & 1.000 \\
1000 & 0.650 & 0.958 & 1.000 & 1.000\\
\hline
\end{tabular*}
\end{table}

The results in Table~\ref{tab1} are based on a single realization of
$\mathbf{Y}_{1000}$ with $\xi=1$ and $\rho=0.01$. We use particle
filters to generate $\{ \mathbf{I}_{1000}^i\dvtx  1 \leq i \leq m \}$
for $m={}$10,000 and compute $\hat\psi_{\mathrm{OR}}$ for $T=200$, 400,
600, 800, 1000, using bootstrap resampling when $\mathrm{cv}_t^2 \geq
c$, for various thresholds $c$. The values of
$\hat\sigma_{\mathrm{OR}}(\hat\psi_{\mathrm{OR}})$ in Table~\ref{tab1}
suggest that variability of $\hat\psi_{\mathrm{OR}}$ is minimized when
$c=2$. Among the 40 simulated values of $\hat\psi_{\mathrm{OR}}$ in
Table~\ref{tab1}, 24 fall within one estimated standard error
$[=\hat\sigma_{\mathrm{OR}} (\hat\psi_{\mathrm{OR}})/\sqrt{m}]$ and 39
within two estimated standard errors of $\psi_T$. This agrees well with
the corresponding numbers 27 and 39, respectively, given by the central
limit theorem in Theorem \ref{thm4a}. Table~\ref{tab2} reports a larger
simulation study involving 500 realizations of $\mathbf{Y}_{1000}$,
with an independent run of the particle filter for each realization.
The table shows confidence intervals for each $T$, with resampling
threshold $c=2$. The results for $\hat\sigma_{\mathrm{OR}}
(\hat\psi_{\mathrm{OR}})$ agree well with the coverage probabilities of
0.683 and 0.954 given by the limiting standard normal distribution in
Theorem~\ref{thm4a}. Table~\ref{tab2} also shows that the
Gilks--Berzuini estimator $\widehat V_T$ described in
Section~\ref{sec5} is very conservative.

\section{\texorpdfstring{Martingales and proof of Theorem \protect\ref{thm3a}}
{Martingales and proof of Theorem 1}}\label{sec3}

In this section we first consider in Section~\ref{sec3.1} the basic
prototype mentioned in the second paragraph of Section~\ref{sec1}, for
which $m(\tilde\psi_T-\psi_T)$ can be expressed as a sum of martingale
differences for a~new class of particle filters. We came up with this
martingale representation, which led to a standard error estimate for
$\tilde\alpha_T$, in \cite{CL11} where we used a particle filter
$\tilde\alpha_T$ to estimate the probability $\alpha$ of a rare event
that $\mathbf{X}_T$ belongs to $\Gamma$. Unlike the present setting of
HMMs, the $\mathbf{X}_T$ is not a vector of hidden states in
\cite{CL11}, where we showed that $m(\tilde\alpha_T-\alpha)$ is a
martingale, thereby proving the unbiasedness of $\tilde\alpha_T$ and
deriving its standard error. In Section~\ref{sec3.1} we extend the
arguments for $m(\tilde\alpha_T-\alpha)$ in \cite {CL11} to
$m(\tilde\psi_T-\psi_T)$ under the assumption that the constants of
proportionality in (\ref{pt}) are easily computable, which we have
called the ``basic prototype'' in Section~\ref{sec1}. We then apply the
results of Section~\ref{sec3.1} to develop in Section~\ref{sec3.2}
martingale approximations for traditional particle filters and use them
to prove Theorem~\ref{thm3a}.

\subsection{Martingale representations for the basic prototype}\label{sec3.1}

We consider here the basic prototype, which requires $T$ to be fixed,
in a more general setting than HMM. Let $(\mathbf{X}_T, \mathbf {Y}_T)$
be a general random vector, with only $\mathbf{Y}_T$ observed, such
that the conditional density $\tilde p_T(\mathbf{X}_T|\mathbf{Y}_T)$ of
$\mathbf{X}_T$ given $\mathbf{Y}_T$ is readily computable. As in
Section~\ref{sec2}, since $\mathbf{Y}_T$ represents the observed data,
we can treat it as constant and simply denote $\tilde
p_T(\mathbf{x}_T|\mathbf{Y}_T)$ by $\tilde p_T(\mathbf{x}_T)$. The
bootstrap resampling scheme in Section~\ref{sec2.1} can be applied to
this general setting and we can use general weight functions
$w_t(\mathbf{x}_t) > 0$, of which (\ref{wt}) is a special case for the
HMM (\ref{hmm}). Because $\tilde p_T(\mathbf{x}_T)$ is readily
computable, so is the likelihood ratio $L_T(\mathbf{x}_T) = \tilde
p_T(\mathbf{x}_T)/\prod_{t=1}^T q_t(x_t|\mathbf{x}_{t-1})$.

This likelihood ratio plays a fundamental role in the theory of the
particle filter estimate $\tilde\alpha_T$ of the rare event probability
$\alpha$ in \cite{CL11}. An unbiased estimate of $\psi_T$ under the
basic prototype is
%
\begin{equation}
\label{23} \tilde\psi_T = m^{-1} \sum
_{i=1}^m L_T \bigl(\widetilde{\mathbf{X}}{}_T^{i} \bigr) \psi \bigl(\widetilde{\mathbf{X}}{}_T^{i} \bigr) H_{T-1}^{i},
\end{equation}
where $H_t^i$ is defined in (\ref{24a}). The unbiasedness follows from
the fact that $m(\tilde\psi_T-\psi_T)$ can be expressed as a sum of
martingale differences. There are in fact two martingale
representations of $m(\tilde\psi_T-\psi_T)$. One is closely related to
that of Del Moral (Chapter~9 of \cite{Del04}); see (\ref{martingale2})
below.\vspace*{1pt} The other is related to the variance estimate $\hat\sigma{}^2
(\hat\psi_T)$ in Theorem \ref{thm3a} and is given in Lemma \ref{lem*}.
Recall the meaning of the notation $E_q$ introduced in the second
paragraph of Section~\ref{sec2.1}. Let $\#_k^{i}$ denote the number of
copies of $\widetilde{\mathbf{X}}{}_k^{i}$ generated from $\{
\widetilde{\mathbf{X}}{}_k^{1},\ldots,\widetilde{\mathbf{X}}{}_k^{m}
\}$ to form the $m$ particles in the $k$th generation. Then,
conditionally, $(\#_k^1,\ldots,\#_k^m) \sim
\operatorname{Multinomial} (m,W_k^1,\ldots, W_k^m)$.

\begin{lem} \label{lem*}
Let $\tilde f_0=\psi_T$ and define for $1 \leq t \leq T$,
%
\begin{equation}
\label{tft} \tilde f_t(\mathbf{x}_t) =
E_q \bigl[\psi(\mathbf{X}_T) L_T(\mathbf
{X}_T)|\mathbf{X}_t=\mathbf{x}_t \bigr].
\end{equation}
Then $\tilde f_t(\mathbf{x}_t) = E_q[\tilde f_T(\mathbf
{X}_T)|\mathbf{X}_t=\mathbf{x}_t]$ and
%
\begin{equation}
\label{29} m(\tilde\psi_T - \psi_T) = \sum
_{j=1}^m \bigl(\varepsilon_1^{j}
+ \cdots+ \varepsilon_{2T-1}^{j} \bigr),
\end{equation}
where
%
\begin{eqnarray}
\label{210} \varepsilon_{2t-1}^{j} & = & \sum
_{i\dvtx A_{t-1}^{i}=j} \bigl[\tilde f_t \bigl(\widetilde{\mathbf{X}}{}_t^{i} \bigr)-\tilde f_{t-1} \bigl(
\mathbf{X}_{t-1}^{i} \bigr) \bigr] H_{t-1}^{i},
\nonumber\\[-9pt]\\[-9pt]
\varepsilon_{2t}^{j} & = & \sum_{i\dvtx A_{t-1}^{i}=j}
\bigl(\#_t^{i}-m W_t^{i} \bigr)
\bigl[\tilde f_t \bigl(\widetilde{\mathbf{X}}{}_t^{i}
\bigr) \widetilde{H}{}_t^{i} -\tilde f_0 \bigr].\nonumber
\end{eqnarray}
Moreover, for each fixed $j$, $\{ \varepsilon_k^{j}, \mathcal{F}_k, 1
\leq k \leq2T-1 \}$ is a martingale difference sequence, where
%
\begin{eqnarray}
\qquad\mathcal{F}_{2t-1} & = & \sigma \bigl( \bigl\{ \widetilde{X}{}_1^{i}\dvtx  1 \leq i \leq m \bigr\}
\cup \bigl\{ \bigl(\mathbf{X}_s^{i},
\widetilde{\mathbf{X}}{}_{s+1}^{i}, A_s^{i}
\bigr)\dvtx  1 \leq s < t, 1 \leq i \leq m \bigr\} \bigr),
\nonumber\\[-8pt]\label{211}  \\[-8pt]
\mathcal{F}_{2t} & = & \sigma \bigl( \mathcal{F}_{2t-1} \cup
\bigl\{ \bigl(\mathbf {X}_t^{i}, A_t^{i}
\bigr)\dvtx  1 \leq i \leq m \bigr\} \bigr)\nonumber
\end{eqnarray}
are the $\sigma$-algebras generated by the random variables associated
with the $m$ particles
just before and just after the
$t$th resampling step, respectively.
\end{lem}

\begin{pf}
Recalling that the ``first generation'' of the $m$ particles consists
of $\widetilde{X}{}_1^{1},\ldots, \widetilde{X}{}_1^{m}$ (before
resampling), $A_t^{i}=j$ if the first component of $\mathbf{X}_t^{i}$
is $\widetilde{X}{}_1^{j}$. Thus, $A_t^{i}$ represents the ``ancestral
origin'' of the ``genealogical particle'' $\mathbf{X}_t^{i}$. It
follows from (\ref{weights}), (\ref{24a}) and simple algebra that for
$1 \leq i \leq m$,
%
\begin{eqnarray}
\label{26} mW_t^{i} & = & H_{t-1}^{i}/
\widetilde{H}{}_t^{i},
\\
\label{27} \sum_{i\dvtx  A_t^{i}=j} \tilde f_t
\bigl(\mathbf{X}_t^{i} \bigr) H_t^{i}
& = & \sum_{i\dvtx
A_{t-1}^{i}=j} \#_t^{i}
\tilde f_t \bigl(\widetilde{\mathbf{X}}{}_t^{i}
\bigr) \widetilde{H}{}_t^{i},
\end{eqnarray}
noting that $\widetilde{\mathbf{X}}{}_t^{i}$ has the same first
component as $\mathbf{X}_{t-1}^{i}$. Multiplying (\ref{26}) by $\tilde
f_t(\widetilde{\mathbf{X}}{}_t^{i}) \widetilde{H}{}_t^i$ and combining
with (\ref{27}) yield
%
\begin{eqnarray}
&& \sum_{t=1}^T  \sum
_{i\dvtx A_{t-1}^{i}=j} \bigl[\tilde f_t \bigl(\widetilde{\mathbf{X}}{}_t^{i} \bigr)- \tilde f_{t-1} \bigl(
\mathbf{X}_{t-1}^{i} \bigr) \bigr] H_{t-1}^{i}\nonumber
\\
&&\quad{}  + \sum_{t=1}^{T-1} \sum
_{i\dvtx A_{t-1}^{i}=j} \bigl(\#_t^{i}-m
W_t^{i} \bigr) \tilde f_t \bigl(\widetilde{\mathbf{X}}{}_t^{i} \bigr) \widetilde{H}{}_t^{i}\label{28}
\\
&&\qquad = \sum_{i\dvtx
A_{T-1}^{i}=j} \tilde f_T \bigl(\widetilde{\mathbf{X}}{}_T^{i} \bigr) H_{T-1}^{i} -
\tilde f_0,\nonumber
\end{eqnarray}
recalling that $A_0^{i}=i$ and $H_0^i=1$. Since $\tilde f_0=\psi _T$
and $\tilde f_T(\mathbf{x}_T)=\psi(\mathbf{x}_T) L_T(\mathbf{x}_T)$ by~(\ref{tft}), it follows from (\ref{23}) that
%
\begin{equation}
\label{psiB} m(\tilde\psi_T-\psi_T) = \sum
_{i=1}^m \tilde f_T \bigl(\widetilde{\mathbf{X}}{}_T^{i} \bigr) H_{T-1}^{i}- m
\tilde f_0.
\end{equation}
Since $\sum_{i=1}^m \#_t^{i} = \sum_{i=1}^m m W_t^{i}=m$ for $1 \leq t
\leq T-1$, combining (\ref{psiB}) with (\ref{28}) yields (\ref{29}).

Let $E_m$ denote expectation under the probability measure induced by
the random variables in the
SMC algorithm with $m$ particles. Then
\[
E_m \bigl(\varepsilon_2^j|\mathcal{F}_1
\bigr)= E_m \bigl\{ \bigl(\#_1^j- m
W_1^j \bigr)| \widetilde{X}{}_1^{1},
\ldots, \widetilde{X}{}_1^{m} \bigr\} \bigl[\tilde
f_1 \bigl(\widetilde{X}{}_1^j \bigr) \widetilde{H}{}_1^j-\tilde f_0 \bigr] = 0
\]
since $E_m (\#_1^j|\mathcal{F}_1)=m W_1^j$. More generally, the
conditional distribution of
$(\mathbf{X}_t^{1},\ldots,\mathbf{X}_t^{m})$ given $\mathcal
{F}_{2t-1}$ is that of $m$ i.i.d. random vectors which take the value
$\widetilde{\mathbf{X}}{}_t^j$ with probability $W_t^j$. Moreover, the
conditional distribution of
$(\widetilde{X}{}_t^{1},\ldots,\widetilde{X}{}_t^{m})$ given $\mathcal
{F}_{2(t-1)}$ is that of $m$ independent random variables such that
$\widetilde{X}{}_t^{i}$ has density function
$q_t(\cdot|\mathbf{X}_{t-1}^{i})$ and, therefore, $E_m[\tilde
f_t(\widetilde{\mathbf{X}}{}_t^{i})|\mathcal{F}_{2(t-1)}] = \tilde
f_{t-1}(\mathbf{X}_{t-1}^{i})$ by (\ref{tft}) and the tower property of
conditional expectations. Hence, in particular,
\[
E_m \bigl(\varepsilon_3^{j}|
\mathcal{F}_2 \bigr) = \sum_{i\dvtx A_1^{i}=j} \bigl\{
E_m \bigl[\tilde f_2 \bigl(\widetilde{\mathbf{X}}{}_2^{i} \bigr)|X_1^{i}
\bigr]- \tilde f_1 \bigl(X_1^{i} \bigr) \bigr\}
H_1^{i} = 0.
\]
Proceeding inductively in this way shows that $\{ \varepsilon_k^{j},
\mathcal{F}_k, 1 \leq k \leq2T-1 \}$ is a martingale difference
sequence for all $j$.
\end{pf}

Without tracing their ancestral origins as in (\ref{210}), we can also
use the successive generations of the $m$ particles to form martingale
differences directly and thereby obtain another martingale
representation of $m(\tilde\psi_T-\psi_T)$. Specifically, the preceding
argument also shows that $\{ (Z_k^1, \ldots, Z_k^m), \mathcal{F}_k, 1
\leq k \leq2T-1 \}$ is a martingale difference sequence, where
%
\begin{eqnarray}
 Z_{2t-1}^{i} & = & \bigl[\tilde
f_t \bigl(\widetilde{\mathbf{X}}{}_t^{i} \bigr)-
\tilde f_{t-1} \bigl(\mathbf{X}_{t-1}^{i} \bigr)
\bigr] H_{t-1}^{i},
\nonumber\\[-8pt]\label{list} \\[-8pt]
Z_{2t}^{i} & = & \tilde f_t \bigl(
\mathbf{X}_t^{i} \bigr) H_t^i-\sum
_{j=1}^m W_t^{j}
\tilde f_t \bigl(\widetilde{\mathbf{X}}{}_t^{j}
\bigr) \widetilde{H}{}_t^{j}.
\nonumber
\end{eqnarray}
Moreover, $Z_k^1, \ldots, Z_k^m$ are conditionally independent given
$\mathcal{F}_{k-1}$. It follows from (\ref{23}), (\ref{26}),
(\ref{list}) and an argument similar to (\ref{28}) that
%
\begin{equation}
\label{martingale2} m(\tilde\psi_T-\psi_T) = \sum
_{k=1}^{2T-1} \bigl(Z_k^1 +
\cdots+ Z_k^m \bigr).
\end{equation}
This martingale representation yields the limiting normal distribution
of $\sqrt{m}(\tilde\psi_T- \psi_T)$ for the basic prototype in the
following.

\begin{cor} \label{thm2} Let
$\sigma^2_{\mathrm{C}} = \tilde\sigma_1^2 + \cdots+
\tilde\sigma_{2T-1}^2$, where
%
\begin{equation}
\label{sigk} \qquad\tilde\sigma_k^2 = %
\cases{E_q \bigl\{ \bigl[\tilde f_t^2(
\mathbf{X}_t)-\tilde f_{t-1}^2(\mathbf
{X}_{t-1}) \bigr] h_{t-1}^*(\mathbf{X}_{t-1}) \bigr\},
&\quad if $k=2t-1$,
\vspace*{2pt}\cr
E_q \bigl\{ \bigl[\tilde f_t(
\mathbf{X}_t) h_t^*(\mathbf{X}_t) - \tilde
f_0 \bigr]^2/h_t^*(\mathbf{X}_t)
\bigr\}, &\quad if $k=2t$.}
\end{equation}
Assume that $\sigma^2_{\mathrm{C}} < \infty$. Then as $m \rightarrow
\infty$,
%
\begin{equation}
\label{clt} \sqrt{m} (\tilde\psi_T-\psi_T)
\Rightarrow N \bigl(0,\sigma^2_{\mathrm{C}} \bigr).
\end{equation}
\end{cor}

The proof of Corollary \ref{thm2} is given in the \hyperref[app]{Appendix}. The main
result of this section is consistent estimation of
$\sigma^2_{\mathrm{C}}$ in the basic prototype. Define for every real
number $\mu$,
%
\begin{eqnarray}
\tilde\sigma{}^2(\mu) & = & m^{-1} \sum
_{j=1}^m \Biggl\{ \sum
_{i\dvtx A_{T-1}^{i}=j} L_T \bigl(\widetilde{\mathbf{X}}{}_T^{i}
\bigr) \psi \bigl(\widetilde{\mathbf{X}}{}_T^{i} \bigr)
H_{T-1}^i
\nonumber\\[-8pt]\label{hsiga} \\[-8pt]
&&\hspace*{42pt}{}- \Biggl[ 1 + \sum_{t=1}^{T-1} \sum
_{i\dvtx A_{t-1}^{i}=j} \bigl(\#_t^{i}-m
W_t^{i} \bigr) \Biggr] \mu \Biggr\}^2,
\nonumber
\end{eqnarray}
and note that by (\ref{210}) and (\ref{28}),
\[
\tilde\sigma{}^2(\psi_T) = m^{-1}
\sum_{j=1}^m \bigl(\varepsilon_1^{j}+
\cdots+\varepsilon_{2T-1}^{j} \bigr)^2.
\]
We next show that $\tilde\sigma{}^2(\psi_T) \stackrel {p}{\rightarrow}
\sigma^2_{\mathrm{C}}$ by making use of the following two lemmas which
are proved in the \hyperref[app]{Appendix}. Since for every fixed $T$ as $m
\rightarrow\infty$, $\tilde\psi_T = \psi_T+O_p(1/\sqrt{m})$ by
(\ref{clt}), it then follows that $\tilde\sigma{}^2(\tilde\psi_T)$ is
also consistent for~$\sigma^2_{\mathrm{C}}$.

\begin{lem} \label{lem1}
Let $1 \leq t \leq T$ and let $G$ be a measurable real-valued function
on the state space. Define $h_t^*$, $\eta_t$ and $\zeta_t$ as in
(\ref{etah}).
\begin{longlist}[(ii)]
\item[(i)] If $E_q[|G(\mathbf{X}_t)|/h_{t-1}^*(\mathbf{X}_{t-1})] <
    \infty$, then as $m \rightarrow\infty$,
%
\begin{equation}
\label{lem11} m^{-1} \sum_{i=1}^m
G \bigl(\widetilde{\mathbf{X}}{}_t^{i} \bigr) \stackrel {p}
{ \rightarrow} E_q \bigl[G(\mathbf{X}_t)/h_{t-1}^*(
\mathbf{X}_{t-1}) \bigr].
\end{equation}

\item[(ii)] If $E_q [|G(\mathbf{X}_t)|/h_t^*(\mathbf{X}_t)] < \infty $,
    then as $m \rightarrow\infty$,
%
\begin{equation}
\label{lem13} m^{-1} \sum_{i=1}^m
G \bigl(\mathbf{X}_t^{i} \bigr) \stackrel{p} {\rightarrow}
E_q \bigl[G(\mathbf {X}_t)/h_t^*(
\mathbf{X}_t) \bigr].
\end{equation}
In particular, for the special case $G=w_t$, (\ref{lem11}) yields $\bar
w_t \stackrel{p}{\rightarrow} \zeta_t^{-1}$ and, hence,
%
\begin{equation}
\label{lem1H} \frac{\widetilde{H}{}_t^i}{h_t^*(\widetilde{\mathbf{X}}{}_t^i)} = \frac
{H_t^i}{h_t^*(\mathbf{X}_t^i)} = \eta_t^{-1}
\bar w_1 \cdots \bar w_t \stackrel{p} {\rightarrow} 1\qquad\mbox{as } m \rightarrow \infty.
\end{equation}
Moreover, if $E_q[|G(\mathbf{X}_t)|/h_{t-1}^*(\mathbf{X}_{t-1})] <
\infty$, then applying (\ref{lem11}) to $|G(\cdot)|\* \mathbf{1}_{\{
|G(\cdot)| > M \}}$ and letting $M \rightarrow\infty$ yield
%
\begin{equation}
\label{lem14} \qquad m^{-1} \sum_{i=1}^m
\bigl|G \bigl(\widetilde{\mathbf{X}}{}_t^{i} \bigr)\bigr|
\mathbf{1}_{\{ |G(\widetilde{\mathbf{X}}{}_t^{i})| > \varepsilon\sqrt{m} \}} \stackrel{p} {\rightarrow} 0\qquad\mbox{as } m \rightarrow
\infty\mbox{ for every } \varepsilon> 0.
\end{equation}
\end{longlist}
\end{lem}

\begin{lem} \label{lem1b}
Let $1 \leq t \leq T$ and let $G$ be a measurable nonnegative valued
function on the state space.
\begin{longlist}[(ii)]
\item[(i)] If $E_q[G(\mathbf{X}_t)/h_{t-1}^*(\mathbf{X}_{t-1})] <
    \infty$, then as $m \rightarrow\infty$,
%
\begin{equation}
\label{lem21} m^{-1} \max_{1 \leq j \leq m} \sum
_{i\dvtx A_{t-1}^{i}=j} G \bigl(\widetilde{\mathbf{X}}{}_t^{i}
\bigr) \stackrel {p} {\rightarrow} 0.
\end{equation}

\item[(ii)] If $E_q[G(\mathbf{X}_t)/h_t^*(\mathbf{X}_t)] <\infty$, then
    as $m \rightarrow\infty$,
%
\begin{equation}
\label{lem22} m^{-1} \max_{1 \leq j \leq m} \sum
_{i\dvtx A_t^{i}=j} G \bigl(\mathbf{X}_t^{i} \bigr)
\stackrel{p} {\rightarrow} 0.
\end{equation}
\end{longlist}
\end{lem}

\begin{cor} \label{thm3}
Suppose $\sigma^2_{\mathrm{C}} < \infty$. Then $\tilde\sigma{}^2(\psi_T)
\stackrel {p}{\rightarrow} \sigma^2_{\mathrm{C}}$ as $m
\rightarrow\infty$.
\end{cor}

\begin{pf}
By (\ref{210}) and (\ref{28}),
\[
\sum_{k=1}^{2T-1}
\varepsilon_k^{j} = \sum_{i\dvtx A_{T-1}^{i}=j}
\tilde f_T \bigl(\widetilde{\mathbf{X}}{}_T^{i}
\bigr) H_{T-1}^i - \Biggl[ 1+ \sum
_{t=1}^{T-1} \sum_{i\dvtx A_{t-1}^{(i)}=j}
\bigl(\#_t^{i}-mW_t^{i} \bigr)
\Biggr] \tilde f_0. %
\]
We make use of Lemmas \ref{lem1} and \ref{lem1b}, (\ref{sigk}) and
mathematical induction to show that for all $1 \leq k \leq2T-1$,
%
\begin{equation}
\label{226} m^{-1} \sum_{j=1}^m
\bigl(\varepsilon_1^{j} + \cdots+ \varepsilon_k^{j}
\bigr)^2 \stackrel{p} {\rightarrow} \sum
_{\ell=1}^k \tilde\sigma_\ell^2.
\end{equation}

By the weak law of large numbers, (\ref{226}) holds for $k=1$. Next
assume that (\ref{226}) holds for $k=2t$ and consider the expansion
%
\begin{eqnarray}
\label{227} && m^{-1} \sum_{j=1}^m
\bigl[ \bigl(\varepsilon_1^{j}+\cdots+\varepsilon_{k+1}^{j}
\bigr)^2- \bigl(\varepsilon _1^{j}+\cdots+
\varepsilon_k^{j} \bigr)^2 \bigr]
\nonumber\\[-8pt]\\[-8pt]
&&\qquad =  m^{-1} \sum_{j=1}^m
\bigl(\varepsilon_{k+1}^{j} \bigr)^2 +
2m^{-1} \sum_{j=1}^m \bigl(
\varepsilon_1^{j} + \cdots+ \varepsilon_k^{j}
\bigr) \varepsilon_{k+1}^{j}.
\nonumber
\end{eqnarray}
Let $C_t^{j} = \{ (i,\ell)\dvtx  i \neq\ell$ and $A_t^{i}=A_t^{\ell}=j
\}$. Suppressing the subscript $t$, let $U_i=[\tilde
f_{t+1}(\widetilde{\mathbf{X}}{}_{t+1}^{i})-\tilde
f_t(\mathbf{X}_t^{i})] H_t^i$. Then $\varepsilon_{2t+1}^{j} =
\sum_{i\dvtx A_t^{i}=j} U_i$ and
%
\begin{equation}
\label{227a} m^{-1} \sum_{j=1}^m
\bigl(\varepsilon_{2t+1}^{j} \bigr)^2 =
m^{-1} \sum_{i=1}^m
U_i^2 + m^{-1} \sum
_{j=1}^m \sum_{(i,\ell) \in
C_t^{j}}
U_i U_\ell.
\end{equation}
By (\ref{sigk}), (\ref{lem1H}) and Lemma \ref{lem1}(i),
%
\begin{equation}
\label{228} m^{-1} \sum_{i=1}^m
U_i^2 \stackrel{p} {\rightarrow} \tilde
\sigma_{2t+1}^2.
\end{equation}
Since $U_1, \ldots, U_m$ are independent mean zero random variables
conditioned on $\mathcal{F}_{2t}$, it follows from (\ref{lem1H}),
Lemmas \ref{lem1}(ii) and \ref{lem1b}(ii) that
%
\begin{eqnarray}
\label{Appendix1}
\qquad && m^{-2} \operatorname{Var}_m \Biggl(
\sum_{j=1}^m \sum_{(i,\ell)\in C_t^{j}} U_i U_\ell \bigg|\mathcal{F}_{2t} \Biggr)\nonumber
\\
&&\qquad=  m^{-2} \sum_{j=1}^m \sum
_{(i,\ell) \in C_t^{j}} E_m \bigl(U_i^2
U_\ell^2| \mathcal{F}_{2t} \bigr)
 \leq
m^{-2} \sum_{j=1}^m \biggl\{
\sum_{i\dvtx A_t^{i}=j} E_m \bigl( U_i^2
| \mathcal{F}_{2t} \bigr) \biggr\}^2
\\
&&\qquad \leq
\Biggl[ m^{-1} \sum_{i=1}^m
E_m \bigl(U_i^2|\mathcal{F}_{2t}
\bigr) \Biggr] \biggl[ m^{-1} \max_{1 \leq j \leq m} \sum
_{i\dvtx A_t^{i}=j} E_m \bigl( U_i^2 |
\mathcal{F}_{2t} \bigr) \biggr] \stackrel{p} {\rightarrow} 0.
\nonumber
\end{eqnarray}
By (\ref{227a})--(\ref{Appendix1}),
%
\begin{equation}
\label{limp} m^{-1} \sum_{j=1}^m
\bigl(\varepsilon_{2t+1}^{j} \bigr)^2 \stackrel{p} {
\rightarrow } \tilde\sigma_{2t+1}^2.
\end{equation}
We next show that
%
\begin{equation}
\label{Appendix2} m^{-1} \sum_{j=1}^m
\bigl(\varepsilon_1^{j} + \cdots+ \varepsilon_{2t}^{j}
\bigr) \varepsilon_{2t+1}^{j} \stackrel{p} {\rightarrow} 0.
\end{equation}
Since $\varepsilon_1^j,\ldots,\varepsilon_{2t}^{j}$ are measurable with
respect to $\mathcal{F}_{2t}$ for $1 \leq j \leq m$ and
$\varepsilon_{2t+1}^{1}, \ldots,\varepsilon_{2t+1}^{m}$ independent
conditioned on $\mathcal{F}_{2t}$, by the induction hypothesis and by~(\ref{lem1H}) and Lemma \ref{lem1b}(ii), it follows that
\begin{eqnarray*}
& & m^{-2} \sum_{j=1}^m
\operatorname{Var}_m \bigl( \bigl(\varepsilon _1^{j}+
\cdots+ \varepsilon_{2t}^{j} \bigr) \varepsilon_{2t+1}^{j}
| \mathcal{F}_{2t} \bigr)
\\
&&\qquad =  m^{-2} \sum
_{j=1}^m \biggl\{ \bigl(\varepsilon_1^{j}
+ \cdots+ \varepsilon _{2t}^{j} \bigr)^2 \sum
_{i\dvtx A_t^{i}=j} E_m \bigl( U_i^2
| \mathcal{F}_{2t} \bigr) \biggr\}
\\
&&\qquad\leq \Biggl[ m^{-1}
\sum_{j=1}^m \bigl(\varepsilon_1^{j}
+ \cdots+ \varepsilon_{2t}^{j} \bigr)^2 \Biggr]
\biggl[ m^{-1} \max_{1 \leq j \leq m} \sum
_{i\dvtx A_t^{i}=j} E_m \bigl(U_i^2|
\mathcal{F}_{2t} \bigr) \biggr] \stackrel{p} {\rightarrow} 0,
\end{eqnarray*}
and therefore (\ref{Appendix2}) holds. By (\ref{227}), (\ref{limp}),
(\ref{Appendix2}) and the induction hypothesis, (\ref{226}) holds for
$k=2t+1$.

Next assume (\ref{226}) holds for $k=2t-1$. Suppressing the subscript
$t$, let $S_i = (\#_t^{i}-mW_t^{i})[\tilde
f_t(\widetilde{\mathbf{X}}{}_t^{i}) \widetilde{H}{}_t^i-\tilde f_0]$.
Then
%
\begin{equation}
\label{235} m^{-1} \sum_{j=1}^m
\bigl(\varepsilon_{2t}^{j} \bigr)^2 =
m^{-1} \sum_{i=1}^m
S_i^2 + m^{-1} \sum
_{j=1}^m \sum_{(i,\ell) \in C_{t-1}^{j}}S_i
S_\ell.
\end{equation}
From Lemma \ref{lem1}(i), (\ref{sigk}) and (\ref{lem1H}), it follows
that
%
\begin{eqnarray}
\label{condV}
&& m^{-1} \sum_{i=1}^m
E_m \bigl(S_i^2|\mathcal{F}_{2t-1}
\bigr)\nonumber
\\
&&\qquad  = m^{-1} \sum_{i=1}^m
\operatorname{Var}_m \bigl(\#_t^i|
\mathcal{F}_{2t-1} \bigr) \bigl[\tilde f_t \bigl(\widetilde{\mathbf{X}}{}_t^i \bigr) \widetilde{H}{}_t^i
- \tilde f_0 \bigr]^2
\\
&&\qquad= m^{-1} \sum_{i=1}^m
\frac{w_t(\widetilde{\mathbf{X}}{}_t^{i})}{\bar w_t} \biggl( 1- \frac{w_t(\widetilde{\mathbf{X}}{}_t^{i})}{m \bar w_t} \biggr) \bigl[\tilde
f_t \bigl(\widetilde{\mathbf{X}}{}_t^{i} \bigr)
\widetilde{H}{}_t^i-\tilde f_0
\bigr]^2 \stackrel{p} {\rightarrow} \tilde\sigma{}_{2t}^2.\nonumber
\end{eqnarray}
Technical arguments given in the \hyperref[app]{Appendix} show that
%
\begin{eqnarray}
\label{eq1} m^{-2} \operatorname{Var}_m \Biggl( \sum
_{i=1}^m S_i^2 \bigg|
\mathcal {F}_{2t-1} \Biggr) & \stackrel{p} {\rightarrow} & 0,
\\
\label{eq2} m^{-2} E_m \Biggl[ \Biggl( \sum
_{j=1}^m \sum_{(i,\ell) \in
C_{t-1}^{j}}
S_i S_\ell \Biggr)^2 \bigg| \mathcal{F}_{2t-1}
\Biggr] & \stackrel{p} {\rightarrow} & 0
\end{eqnarray}
and 
\begin{equation}\label{appendix4} m^{-2} E_m \Biggl[ \Biggl\{ \sum
_{j=1}^m \bigl(\varepsilon_1^{j}
+ \cdots+ \varepsilon_{2t-1}^{j} \bigr) \varepsilon_{2t}^{j}
\Biggr\}^2 \bigg| \mathcal {F}_{2t-1} \Biggr]  \stackrel{p} {
\rightarrow}  0
\end{equation}
under (\ref{226}) for $k=2t-1$. By (\ref{227}) and
(\ref{235})--(\ref{appendix4}), (\ref{226}) holds for $k=2t$. The
induction proof is complete and Corollary \ref{thm3} holds.
\end{pf}

Let $\Psi(\mathbf{x}_T) = \psi(\mathbf{x}_T)-\psi_T$ and define
$\widetilde\Psi_T$ as in (\ref{23}) with $\Psi$ in place of $\psi$.
Then replacing $\tilde\psi_T$ by $\widetilde\Psi_T$ in Corollaries
\ref{thm2} and \ref{thm3} shows that as $m \rightarrow\infty$,
%
\begin{equation}
\label{tPsi} \sqrt{m} \widetilde\Psi_T \Rightarrow N \bigl(0,
\sigma^2 \bigr), \qquad \tilde\sigma_*^2(0) \stackrel{p} {
\rightarrow} \sigma^2,
\end{equation}
where $\tilde\sigma_*^2(\mu)$ is defined by (\ref{hsiga}) with $\psi$
replaced by $\Psi$, noting that replacing $\psi$ by $\Psi$ in the
definition of $\sigma^2_{\mathrm{C}}$ in Corollary \ref{thm2} gives
$\sigma^2$ defined in (\ref{tsigk}) as $f_0=0$.

\subsection{\texorpdfstring{Approximation of $\hat\psi_T-\psi_T$ by $\widetilde
\Psi_T$ and proof of Theorem \protect\ref{thm3a}}
{Approximation of psi T - psi T by Psi T and proof of Theorem 1}}\label{sec3.2}

We now return to the setting of Section~\ref{sec2.1}, in which we
consider the HMM (\ref{hmm}) and use weight functions of the form
(\ref{wt}) and the particle filter (\ref{23a}) in lieu of (\ref{23}).
We make use of the following lemma to extend (\ref{tPsi}) to
$\hat\psi_T-\psi_T$.

\begin{lem} \label{lem31}
Let $w_t$ be of the form (\ref{wt}) for $1 \leq t \leq T$ and define
$\eta_T$ and $\zeta_T$ by (\ref{etah}) and $\widetilde\Psi _T$ by
(\ref{23}) with $\psi$ replaced by $\Psi= \psi-\psi_T$. Then
%
\begin{eqnarray}
\label{l21} \tilde p_T(\mathbf{x}_T)\, \bigl[\! &=&
\tilde p_T(\mathbf {x}_T|\mathbf{Y}_T) \bigr]
= \eta_T^{-1} \prod_{t=1}^T
\bigl[p_t(x_t|x_{t-1}) g_t(Y_t|x_t)
\bigr],
\\
\label{l22} \qquad L_T(\mathbf{x}_T) & = &
\eta_T^{-1} \prod_{t=1}^T
w_t(\mathbf {x}_t),
\\
\label{l23} \hat\psi_T - \psi_T & = & (m \bar
w_T)^{-1} \sum_{i=1}^m
\Psi \bigl(\widetilde{\mathbf{X}}{}_T^i \bigr)
w_T \bigl(\widetilde{\mathbf{X}}{}_T^i \bigr) =
( \bar w_1 \cdots\bar w_T)^{-1}
\eta_T \widetilde\Psi_T.
\end{eqnarray}
\end{lem}

\begin{pf}
By (\ref{wt}) and (\ref{etah}),
\[
\eta_T = \int\prod_{t=1}^T
\bigl[w_t(\mathbf{x}_t) q_t(x_t|
\mathbf {x}_{t-1}) \bigr] \,d \nu_X(\mathbf{x}_T)
= \int\prod_{t=1}^T \bigl[p_t(x_t|x_{t-1})g_t(Y_t|x_t)
\bigr] \,d \nu_X(\mathbf{x}_T),
\]
and (\ref{l21}) follows from (\ref{pt}). Moreover, (\ref{l22}) follows
from (\ref{wt}), (\ref{kappat}) and (\ref{l21}). The first equality in
(\ref{l23}) follows from (\ref{24a}) and (\ref{23a}). By (\ref{24a})
and~(\ref{l22}),
%
\begin{equation}
\label{LH} L_T \bigl(\widetilde{\mathbf{X}}{}_T^i
\bigr) H_{T-1}^i = \biggl( \frac{\bar w_1
\cdots\bar w_T}{\eta_T} \biggr)
\frac{w_T(\widetilde{\mathbf{X}}{}_T^i)}{\bar w_T}.
\end{equation}
Hence, the second equality in (\ref{l23}) follows from (\ref{23}).
\end{pf}

The following lemma, proved in the \hyperref[app]{Appendix}, is used to prove
Theoerem~\ref{thm3a}. Although it resembles Lemma \ref{lem1b}, its
conclusions are about the square of the sums in Lemma \ref{lem1b} and,
accordingly, it assumes finiteness of second (instead of first)
moments.

\begin{lem} \label{lem7} Let $1 \leq t \leq T$, $G$ be a measurable
function on the state space, and $\Gamma_t(\mathbf{x}_t)$ be defined as
in (\ref{etah}).
\begin{longlist}[(ii)]
\item[(i)] If $E_q[G^2(\mathbf{X}_t) \Gamma_{t-1}(\mathbf{X}_{t-1})] <
\infty$, then
\[
m^{-1} \sum_{j=1}^m \biggl[
\sum_{i\dvtx A_{t-1}^i=j} G \bigl(\widetilde{\mathbf{X}}{}_t^i
\bigr) \biggr]^2 =O_p(1)\qquad\mbox{as } m \rightarrow
\infty.
\]

\item[(ii)] If $E_q[G^2(\mathbf{X}_t) \Gamma_t(\mathbf{X}_t)] < \infty$, then
\[
m^{-1} \sum_{j=1}^m \biggl[
\sum_{i\dvtx A_t^i=j} G \bigl(\mathbf{X}_t^i
\bigr) \biggr]^2 = O_p(1)\qquad\mbox{as } m \rightarrow
\infty.
\]
\end{longlist}
\end{lem}

\begin{pf*}{Proof of Theorem \ref{thm3a}}
By (\ref{lem1H}), (\ref{tPsi}) and (\ref{l23}),
%
\begin{equation}
\label{clthat} \sqrt{m} (\hat\psi_T-\psi_T) =
\bigl(1+o_p(1) \bigr) \sqrt{m} \widetilde \Psi_T\quad\Rightarrow\quad N \bigl(0,\sigma^2 \bigr).
\end{equation}
Since\vspace*{-1pt} $\tilde\sigma_*^2(0) = m^{-1} \sum_{j=1}^m \{ \sum_{i\dvtx
A_{T-1}^i=j} L_T(\widetilde{\mathbf{X}}{}_T^i) H_{T-1}^i \Psi
(\widetilde{\mathbf{X}}{}_T^i) \}^2 \stackrel{p}{\rightarrow} \sigma^2$
by (\ref{hsiga}) and (\ref{tPsi}), and since $\tilde\sigma_*^2(0)=
(1+o_p(1)) \hat\sigma{}^2(\psi_T)$ in view of (\ref{420}),
(\ref{lem1H}) and~(\ref{LH}),
%
\begin{equation}
\label{wh2sig} \hat\sigma{}^2(\psi_T) \stackrel{p} {
\rightarrow} \sigma^2.
\end{equation}
Letting\vspace*{-1pt} $c_j = \sum_{i\dvtx A_{T-1}^i=j}
\frac{w_T(\widetilde{\mathbf{X}}{}_T^i)}{\bar
w_T}[\psi(\widetilde{\mathbf{X}}{}_T^i)-\psi_T]$, $b_j = \sum_{i\dvtx
A_{T-1}^i=j} w_T(\widetilde{\mathbf{X}}{}_T^i)/\bar w_T$ and
$a=\hat\psi_T-\psi_T$, it follows from (\ref{420}) that
%
\begin{equation}
\label{sigexp} \qquad\hat\sigma{}^2(\hat\psi_T) =
\frac{1}{m} \sum_{j=1}^m
(c_j+a b_j)^2 = \frac{1}{m} \sum
_{j=1}^m c_j^2 +
\frac{2a}{m} \sum_{j=1}^m
b_j c_j + \frac{a^2}{m} \sum
_{j=1}^m b_j^2.
\end{equation}
Since $a = O_p(m^{-1/2})$ by (\ref{clthat}), $m^{-1} \sum_{j=1}^m c_j^2
= \hat\sigma{}^2(\psi_T) \stackrel{p}{\rightarrow} \sigma^2$ by
(\ref{wh2sig}), and since $m^{-1} \sum_{j=1}^m b_j^2 = O_p(1)$ by
Lemma~\ref{lem7} (with $G=w_t$),
\[
\Biggl| \frac{2a}{m} \sum_{j=1}^m
b_j c_j \Biggr| \leq|a| \Biggl( m^{-1} \sum
_{j=1}^m b_j^2 +
m^{-1} \sum_{j=1}^m
c_j^2 \Biggr) \stackrel{p} {\rightarrow} 0.
\]
Hence, by (\ref{sigexp}), $\hat\sigma{}^2(\hat\psi_T)
\stackrel{p}{\rightarrow} \sigma^2$.
\end{pf*}

\section{\texorpdfstring{Extensions and proof of Theorem \protect\ref{thm4a}}
{Extensions and proof of Theorem 2}}\label{sec4}

In this section we first extend the results of Section~\ref{sec2.1} to
the case where residual Bernoulli resampling is used in lieu of
bootstrap resampling at every stage. We then consider occasional
bootstrap resampling and prove Theorem \ref{thm4a}, which we also
extend to the case of occasional residual Bernoulli resampling.

\subsection{Residual Bernoulli resampling}\label{sec4.1}

The \textit{residual resampling} scheme introduced in
\cite{Bak85,Bak87} often leads to smaller asymptotic variance than that
of bootstrap resampling; see \cite{Cho04,LC98}. We consider here the
residual Bernoulli scheme given in \cite{CDL99}, which has been shown
in Section~2.4 of \cite{DM08} to yield a consistent and asymptotically
normal particle filtering estimate of $\psi_T$. To implement residual
Bernoulli resampling, we modify bootstrap resampling at stage $t$ as
follows: let $M(1)=m$ and let $\xi_t^1, \ldots, \xi_t^{M(t)}$ be
independent Bernoulli random variables conditioned on $(M(t),W_t^i)$
satisfying\vspace*{1pt} $P \{ \xi_t^i = 1 | \mathcal{F}_{2t-1} \} = M(t) W_t^{i} -
\lfloor M(t) W_t^{i} \rfloor$. For each $1 \leq i \leq M(t)$ and $t
\leq T-1$, let $\widetilde{H}{}_t^i = H_{t-1}^i/(M(t) W_t^i)$ and make
$\#_t^{i}:=\lfloor M(t) W_t^i \rfloor+ \xi_t^i$ copies of
$(\widetilde{\mathbf{X}}{}_t^{i}, A_{t-1}^{i}, \widetilde{H}{}_t^i)$.
These\vspace*{-1pt} copies constitute an augmented sample $\{ (\mathbf{X}_t^{j},
A_t^{j}, H_t^{j})\dvtx\break  1 \leq j \leq M(t+1) \}$, where $M(t+1) =
\sum_{i=1}^{M(t)} \#_t^{i}$. Estimate $\psi_T$ by
\[
\hat\psi_{T,\mathrm{R}}:= \bigl[M(T) \bar w_T
\bigr]^{-1} \sum_{i=1}^{M(T)} \psi
\bigl(\widetilde{\mathbf{X}}{}_T^i \bigr) w_T
\bigl(\widetilde{\mathbf{X}}{}_T^i \bigr).
\]
Let $\widetilde\Psi_{T,\mathrm{R}}$ be (\ref{23}) with $\Psi= \psi-
\psi_T$ in place of $\psi$ and $M(T)$ in place of $m$. Since
$H_{t-1}^i$ is still given by (\ref{24a}) and Lemma \ref{lem1} still
holds for residual Bernoulli resampling, it follows from (\ref{lem1H})
and (\ref{LH}) that
%
\begin{eqnarray}
M(T) (\hat\psi_{T,{\mathrm{R}}}-\psi_T) &=& \bar
w_T^{-1} \sum_{i=1}^{M(T)}
\Psi \bigl(\widetilde{\mathbf{X}}{}_T^i \bigr)
w_T \bigl(\widetilde{\mathbf{X}}{}_T^i \bigr)\nonumber
\\
& = & \biggl( \frac{\eta_T}{\bar w_1 \cdots\bar w_T} \biggr) \sum_{i=1}^{M(T)}
\Psi \bigl(\widetilde{\mathbf{X}}{}_T^i \bigr)
L_T \bigl(\widetilde{\mathbf{X}}{}_T^i \bigr)
h_{t-1} \bigl(\mathbf{X}_{t-1}^i \bigr)\label{rs1}  
\\
&=& \bigl(1+o_p(1) \bigr) M(T) \widetilde
\Psi_{T, \mathrm{R}}.\nonumber
\end{eqnarray}
Because replacing $\psi$ by $\Psi$ in (\ref{tft}) modifies $\tilde f_t$
to $f_t$ as defined in (\ref{25}), we have for residual Bernoulli
resampling the following analog of (\ref{martingale2}):
%
\begin{equation}
\label{rs2} M(T) \widetilde\Psi_{T, \mathrm{R}} = \sum
_{k=1}^{2T-1} \bigl(Z_{k,\mathrm{R}}^1 +
\cdots+ Z_{k,\mathrm{R}}^{M(\lceil
k/2 \rceil)} \bigr),
\end{equation}
where $Z_{2t-1,\mathrm{R}}^i =
[f_t(\widetilde{\mathbf{X}}{}_t^i)-f_{t-1}(\mathbf{X}_{t-1}^i)]
H_{t-1}^i$, $Z_{2t,\mathrm{R}}^i = (\#_t^i-M(t) W_t^i)
f_t(\widetilde{\mathbf{X}}{}_t^i) \widetilde{H}{}_t^i$.

Let $\sigma^2_{\mathrm{R}} = \sum_{k=1}^{2T-1}
\sigma_{k,\mathrm{R}}^2$, with
\[
\sigma^2_{2t-1,\mathrm{R}} = \sigma^2_{2t-1},
\qquad \sigma_{2t,\mathrm{R}}^2 = E_q \bigl[ \gamma
\bigl(w_t^*(\mathbf{X}_t) \bigr)f_t^2(
\mathbf{X}_t) h_t^*(\mathbf{X}_t) \bigr],
\]
in which $\gamma(x)=(x-\lfloor x \rfloor)(1-x+\lfloor x \rfloor)/x$,
$w_t^*(\mathbf{x}_t) = h_{t-1}^*(\mathbf{x}_{t-1})/h_t^*(\mathbf {x}_t)
=\zeta_t w_t(\mathbf{x}_t)$ and $\sigma^2_{2t-1}$ is defined in Theorem
\ref{thm3a}. Since $\operatorname{Var}(\xi) = \gamma(x)x$ when $\xi$ is
a~Bernoulli($x-\lfloor x \rfloor$) random variable,
\[
\operatorname{Var} \Biggl( \sum_{i=1}^{M(t)}
Z_{2t,{\mathrm{R}}}^i \bigg| \mathcal{F}_{2t-1} \Biggr) = \sum
_{i=1}^{M(t)} \gamma \bigl(M(t)
W_t^i \bigr) M(t) W_t^i
f_t^2 \bigl(\widetilde{\mathbf{X}}{}_t^i
\bigr) \bigl(\widetilde{H}{}_t^i \bigr)^2
\]
and, therefore, it follows from Lemma \ref{lem1} that
%
\begin{equation}
\label{rs3} m^{-1} \operatorname{Var} \Biggl( \sum
_{i=1}^{M(\lceil k/2 \rceil)} Z_{k,\mathrm{R}}^i \bigg|
\mathcal{F}_{k-1} \Biggr) \stackrel {p} {\rightarrow}
\sigma^2_{k,\mathrm{R}}, \qquad 1 \leq k \leq2T-1.
\end{equation}
By (\ref{rs1})--(\ref{rs3}) and a similar modification of (\ref{29})
for residual Bernoulli resampling, Theorem \ref{thm3a} still holds with
$\hat\psi_T$ replaced by $\hat\psi_{T,{\mathrm{R}}}$; specifically,
%
\begin{equation}
\label{rs4} \sqrt{m}(\hat\psi_{T, \mathrm{R}}-\psi_T) \Rightarrow
N \bigl(0,\sigma ^2_{\mathrm{R}} \bigr), \qquad\hat
\sigma^2(\hat\psi_{T, \mathrm{R}}) \stackrel {p} {\rightarrow}
\sigma^2_{\mathrm{R}},
\end{equation}
where $\hat\sigma{}^2(\mu)$ is defined in (\ref{420}). The following
example shows the advantage of residual resampling over bootstrap
resampling.

\begin{exa*}
Consider the bearings-only tracking problem in \cite{GSS93}, in which
$X_t=(X_{t1},\ldots, X_{t4})'$ and $(X_{t1},X_{t3})$ represents the
coordinates of a ship and $(X_{t2}, X_{t4})$ the velocity at time $t$.
A person standing at the origin observes that the ship is at angle
$Y_t$ relative to him. Taking into account measurement errors and
random disturbances in velocity of the ship, and letting $\Phi$ be a $4
\times4$ matrix with $\Phi_{12} = \Phi_{34} =1 = \Phi_{ii}$ for $1 \leq
i \leq4$ and all other entries 0, $\Gamma $ be a $4 \times2$ matrix
with $\Gamma_{11} = \Gamma_{32} = 0.5$, $\Gamma_{21} = \Gamma_{42} = 1$
and all other entries 0, the dynamics can be described by state-space
model
\[
X_t = \Phi X_{t-1} + \Gamma z_t, \qquad
Y_t = \tan^{-1}(X_{t3}/X_{t1})+u_t,
\]
where $X_{11} \sim N (0,0.5^2)$,  $X_{12} \sim N(0,0.005^2)$, $X_{13}
\sim N(0.4,0.3^2)$,  $X_{14} \sim\break  N(-0.05,0.01^2)$, $z_{t+1} \sim
N(0,0.001^2 \mathbf{I}_2)$ and $u_t \sim N(0,0.005^2)$ are independent
random variables for $t \geq1$. The quantities of interest are
$E(X_{T1}|\mathbf{Y}_T)$ and $E(X_{T3}|\mathbf{Y}_T)$. To implement the
particle filters, we let $q_t=p_t$ for $t \geq2$. Unlike \cite{GSS93}
in which $q_1=p_1$ as well, we use an importance density $q_1$ that
involves $p_1$~and~$Y_1$ to generate the initial location
$(X_{11},X_{13})$. Let $\xi$ and $\zeta$ be independent standard normal
random variables, and $r = \tan(Y_1+0.005 \xi)$. Since $u_1$ has small
variance, we can estimate $X_{13}/X_{11}$ well by $r$. This suggests
choosing $q_1$ to be degenerate bivariate normal with support $y=rx$,
with $(y,x)$ denoting $(X_{13},X_{11})$, so that its density function
on this line is proportional to $\exp\{-(x-\mu)^2/(2 \tau) \}$, where
$\mu =0.4r/(0.36+r^2)$, $\tau=0.09/(0.36+r^2)$. Thus, $q_1$ generates
$X_{11}=\mu+\sqrt{\tau} \zeta$ and $X_{13} = r X_{11}$, but still
follows $p_1$ to generate $(X_{12},X_{14})$. By (\ref{wt}),
\[
w_1(x_1) = |x_{11}| \sqrt{\tau}
\bigl(1+r^2 \bigr) \exp \biggl\{ -\frac{x_{11}^2}{2(0.5)^2}-
\frac{(x_{13}-0.4)^2}{2(0.3)^2} + \frac
{\zeta^2}{2} \biggr\},
\]
noting that $0.005 |x_{11}| \sqrt{\tau} (1+r^2)$ is the Jacobian of the
transformation of $(\xi,\zeta)$ to $(x_{11},x_{13})$. For $t \geq2$,
(\ref{wt}) yields the resampling weights $w_t(\mathbf{x}_t) = \exp\{
-[Y_t-\tan ^{-1}(x_{t3}/x_{t1})]^2/[2(0.005)^2] \}$. We use
$m={}$10,000 particles to estimate $E((X_{T1},X_{T3})|\mathbf {Y}_T)$
by particle filters, using bootstrap (boot) or residual (resid)
Bernoulli resampling at every stage, for different values of $T$.
Table~\ref{tab3}, which reports results based on a single realization
of $\mathbf{Y}_T$, shows that residual Bernoulli resampling has smaller
standard errors than bootstrap resampling. Table~\ref{tab3} also
considers bootstrap resampling that uses $q_t=p_t$ for $t \geq1$, as in
\cite{GSS93} and denoted by $\operatorname{boot}(P)$, again with $m={}$10,000
particles. Although boot\vadjust{\goodbreak} and $\operatorname{boot}(P)$ differ only in the choice of the
importance densities at $t=1$, the standard error of boot is
substantially smaller than that of $\operatorname{boot}(P)$ over the entire range $4
\leq T \leq24$. Unlike Section~\ref{sec2.3}, the standard error
estimates in Table~\ref{tab3} use a~sample splitting refinement that
is described in the first paragraph of Section~\ref{sec5}.
\end{exa*}\vspace*{-3pt} 

\begin{table}
\tabcolsep=0pt
\caption{Simulated values of $\hat\psi_T \pm\hat\sigma (\hat\psi_T)
/\sqrt{m}$ (for bootstrap resampling) and $\hat\psi_{T,\mathrm{R}}
\pm\hat\sigma(\hat\psi_{T,\mathrm{R}})/\sqrt{m}$ (for residual
Bernoulli resampling)}\label{tab3}
\begin{tabular*}{\tablewidth}{@{\extracolsep{\fill}}@{}lcccccc@{}}
\hline
& \multicolumn{3}{c}{$\bolds{E(X_{T1}|\mathbf{Y}_T)}$} & \multicolumn{3}{c}{$\bolds{E(X_{T3}|\mathbf{Y}_T)}$}\\[-6pt]
& \multicolumn{3}{c}{\hrulefill} & \multicolumn{3}{c}{\hrulefill}\\
$\bolds{T}$ & $\bolds{\operatorname{boot}(P)}$ & \textbf{boot} & \textbf{resid} & $\bolds{\operatorname{boot}(P)}$ & \textbf{boot} & \textbf{resid}\\
\hline
\phantom{0}4 & $-$0.515 & $-$0.536 & $-$0.519 & \phantom{$-$}0.091 & \phantom{$-$}0.098 & \phantom{$-$}0.095 \\
& $\pm$0.026 & $\pm$0.009 & $\pm$0.007 & $\pm$0.007 & $\pm$0.002 & $\pm$0.001 \\[2pt]
\phantom{0}8 & $-$0.505 & $-$0.515 & $-$0.497 & $-$0.118 & $-$0.121 & $-$0.117 \\
& $\pm$0.043 & $\pm$0.012 & $\pm$0.011 & $\pm$0.008 & $\pm$0.003 & $\pm$0.003 \\[2pt]
12 & $-$0.503 & $-$0.532 & $-$0.508 & $-$0.327 & $-$0.346 & $-$0.331 \\
& $\pm$0.042 & $\pm$0.014 & $\pm$0.010 & $\pm$0.027 & $\pm$0.009 & $\pm$0.007 \\[2pt]
16 & $-$0.506 & $-$0.538 & $-$0.516 & $-$0.538 & $-$0.572 & $-$0.549 \\
& $\pm$0.044 & $\pm$0.015 & $\pm$0.010 & $\pm$0.046 & $\pm$0.016 & $\pm$0.011 \\[2pt]
20 & $-$0.512 & $-$0.537 & $-$0.532 & $-$0.747 & $-$0.783 & $-$0.777 \\
& $\pm$0.045 & $\pm$0.016 & $\pm$0.011 & $\pm$0.066 & $\pm$0.023 & $\pm$0.017 \\[2pt]
24 & $-$0.512 & $-$0.551 & $-$0.540 & $-$0.950 & $-$1.023 & $-$1.002 \\
& $\pm$0.047 & $\pm$0.016 & $\pm$0.012 & $\pm$0.088 & $\pm$0.029 & $\pm$0.022\\
\hline
\end{tabular*}
\end{table}

\subsection{\texorpdfstring{Proof of Theorem \protect\ref{thm4a}}
{Proof of Theorem 2}}\label{sec4.2}
First consider the following modification of~$\widetilde\Psi_T$:
%
\begin{equation}
\label{423} \qquad \widetilde\Psi_{\mathrm{OR}}:= m^{-1} \sum
_{i=1}^m L_T \bigl(\widetilde{\mathbf{X}}{}_T^i \bigr) \Psi \bigl(\widetilde{\mathbf{X}}{}_T^i \bigr) H_{\tilde\tau(T)}^i =
\eta_T^{-1} \bar v_{\tau_1} \cdots\bar
v_{\tau_{r+1}} (\hat\psi_{\mathrm{OR}}-\psi_T).
\end{equation}
Analogous to (\ref{martingale2}), $m \widetilde\Psi_{\mathrm{OR}} =
\sum_{k=1}^{2r+1} (Z_k^1 + \cdots + Z_k^m)$, where
\begin{eqnarray*}
Z_{2s-1}^i & = & \bigl[f_{\tau_s} \bigl(\widetilde{\mathbf{X}}{}_{\tau_s}^i \bigr)- f_{\tau_{s-1}} \bigl(
\mathbf{X}_{\tau_{s-1}}^i \bigr) \bigr] H_{\tau_{s-1}}^i,
\\[-3pt] 
Z_{2s}^i & = & f_{\tau_s} \bigl(\widetilde{\mathbf{X}}{}_{\tau_s}^i \bigr) \widetilde{H}{}_{\tau_s}^i-
\sum_{j=1}^m V_{\tau_s}^j
f \bigl(\widetilde{\mathbf{X}}{}_{\tau_s}^j \bigr) \widetilde{H}{}_{\tau_s}^j.
\end{eqnarray*}
Let $\tau_0^*=0$. Since $\bar v_{\tau_s} \stackrel{p}{\rightarrow}
\prod_{t=\tau _{s-1}^*+1}^{\tau_s^*} \zeta_t$ by (\ref{rstar}) and
Lemma \ref{lem1}, it follows from (\ref{etah}) and (\ref{423}) that
\[
\sqrt{m}(\hat\psi_{\mathrm{OR}}-\psi_T) = \bigl(1+o_p(1)
\bigr) \sqrt{m} \widetilde\Psi_{\mathrm{OR}}\quad\Rightarrow\quad N \bigl(0,
\sigma_{\mathrm{OR}}^2 \bigr),
\]
similarly to (\ref{clthat}), where
%
\begin{eqnarray}
\sigma^2_{\mathrm{OR}} & = & \sum
_{s=1}^{r^*+1} E_q \bigl\{
\bigl[f_{\tau^*_s}^2(\mathbf{X}_{\tau
^*_s})-f_{\tau^*_{s-1}}^2
(\mathbf{X}_{\tau^*_{s-1}}) \bigr] h^*_{\tau^*_{s-1}}(\mathbf{X}_{\tau
^*_{s-1}})\bigr\}
\nonumber\\[-8pt]\label{sigOR} \\[-8pt]
&&{}+ \sum_{s=1}^{r^*}
E_q \bigl[f_{\tau^*_s}^2(\mathbf{X}_{\tau
^*_s})
h_{\tau^*_s}^* (\mathbf{X}_{\tau_s^*}) \bigr].\nonumber
\end{eqnarray}
Moreover, analogous to (\ref{29}), $m \widetilde\Psi_{\mathrm{OR}} =
\sum_{j=1}^m (\varepsilon_1^j+\cdots+\varepsilon_{2r+1}^j)$, where
\begin{eqnarray*}
\varepsilon_{2s-1}^j &= & \sum
_{i\dvtx  A_{\tau_{s-1}}^i =j} \bigl[f_{\tau
_s} \bigl(\widetilde{\mathbf{X}}{}_{\tau_s}^i \bigr)- f_{\tau_{s-1}} \bigl(
\mathbf{X}_{\tau_{s-1}}^i \bigr) \bigr] H_{\tau_{s-1}}^i,
\\
\varepsilon_{2s}^j &= & \sum_{i\dvtx A_{\tau_{s-1}}^i =j}
\bigl(\#_{\tau_s}^i -mV_{\tau_s}^i \bigr)
f_{\tau_s} \bigl(\widetilde{\mathbf{X}}{}_{\tau_s}^i
\bigr) \widetilde{H}{}_{\tau_s}^i,
\end{eqnarray*}
therefore, the proof of $\hat\sigma_{\mathrm{OR}}^2(\hat\psi
_{\mathrm{OR}}) \stackrel{p}{\rightarrow} \sigma^2_{\mathrm{OR}}$ is
similar to that of Theorem~\ref{thm3a}.

\subsection{\texorpdfstring{Occasional residual resampling and assumption (\protect\ref{rstar})}
{Occasional residual resampling and assumption (2.12)}}\label{sec4.3}
In the case of occasional residual resampling, Theorem \ref{thm4a}
still holds with $\hat\psi_{\mathrm{OR}}$ replaced by
$\hat\psi_{\mathrm{ORR}}$ and with $\sigma ^2_{\mathrm{OR}}$ replaced
by
\begin{eqnarray*}
\sigma_{\mathrm{ORR}}^2 & = & \sum_{s=1}^{r^*+1}
E_q \bigl\{ \bigl[f_{\tau
^*_s}^2(\mathbf{X}_{\tau^*_s})- f_{\tau^*_{s-1}}^2(\mathbf{X}_{\tau^*_{s-1}})
\bigr] h^*_{\tau^*_{s-1}}(\mathbf{X}_{\tau^*_{s-1}}) \bigr\}
\\
&&{}+ \sum_{s=1}^{r^*} E_q \Biggl[
\gamma \Biggl( \prod_{t=\tau
^*_{s-1}+1}^{\tau^*_s}
w_t^*( \mathbf{X}_t) \Biggr) f_{\tau^*_s}^2
( \mathbf{X}_{\tau^*_s}) h_{\tau^*_s}^* (\mathbf{X}_{\tau^*_s})
\Biggr].
\end{eqnarray*}
In particular, $\sqrt{m}(\hat\psi_{\mathrm{ORR}}-\psi_T) \Rightarrow
N(0,\sigma^2_{\mathrm{ORR}})$ and
$\hat\sigma{}^2_{\mathrm{OR}}(\hat\psi_{\mathrm{ORR}}) \stackrel
{p}{\rightarrow} \sigma^2_{\mathrm{ORR}}$.

Assumption (\ref{rstar}) is often satisfied in practice. In particular,
if one follows Liu~\cite{Liu08} and resamples at stage $t$ whenever
%
\begin{equation}
\label{cvt} {\mathrm{cv}}_t^2:= m^{-1} \sum
_{i=1}^m \bigl(m V_t^i
\bigr)^2-1 = \frac{\sum_{i=1}^m
[v_t(\widetilde{\mathbf{X}}{}_t^i)-\bar v_t]^2}{m \bar v_t^2} \geq c,
\end{equation}
where $c>0$ is a prespecified threshold for the coefficient of
variation, then (\ref{rstar}) can be shown to hold by making use of the
following lemma.

\begin{lem} \label{lem0}
Let resampling be carried out at stage $t$ whenever ${\mathrm{cv}}_t^2
\geq c$. Then (\ref{rstar}) holds with
%
\begin{eqnarray}\label{tauq}
\tau^*_s = \inf \biggl\{ t>\tau^*_{s-1}\dvtx
E_q \biggl( \frac{[ \prod_{k=\tau^*_{s-1}+1}^t
w_k^*(\mathbf{X}_k)]^2}{h_{\tau^*_{s-1}}^*
(\mathbf{X}_{\tau^*_{s-1}})} \biggr)-1 \geq c \biggr\}
\nonumber\\[-16pt]\\[-2pt]
\eqntext{\mbox{for } 1 \leq s \leq r^*,}
\end{eqnarray}
where $\tau_0^*=0$ and $r^* = \max\{ s\dvtx  \tau^*_s < T \}$, provided
that
%
\begin{equation}
\label{CVcond} E_q \biggl( \frac{[ \prod_{k=\tau^*_{s-1}+1}^{\tau^*_s}
w_k^*(\mathbf{X}_k)]^2}{h^*_{\tau^*_{s-1}}(\mathbf{X}_{\tau
^*_{s-1}})} \biggr) -1 \neq c\qquad\mbox{for } 1 \leq s \leq r^*.
\end{equation}
\end{lem}

\begin{pf} Let\vspace*{-2pt} $\tau_{s-1}^*+1 \leq\ell\leq\tau_s^*$ and
$G(\mathbf{x}_\ell )=\prod_{k=\tau_{s-1}^*+1}^\ell
w_k^*(\mathbf{x}_k)$. Apply Lemma \ref{lem1}(i) to $G$ and $G^2$, with
$t-1$ in the subscript replaced by the most recent resampling time
$\tau^*_{s-1}$. It then follows from (\ref{cvt}) that
%
\begin{equation}
\label{cvl} {\mathrm{cv}}_\ell^2 \stackrel{p} {\rightarrow}
E_q \bigl[G^2(\mathbf{X}_\ell
)/h_{\tau_{s-1}^*}^* (\mathbf{X}_{\tau_{s-1}^*}) \bigr] -1\qquad\mbox{on } \bigl\{
\tau_{s-1}=\tau _{s-1}^* \bigr\}.
\end{equation}
In\vspace*{-2pt} view of (\ref{tauq}) and (\ref{CVcond}), it follows from (\ref{cvl})
with $\ell=\tau_s^*$ that $P_m \{ \tau_s > \tau_s^*,
\tau_{s-1}=\tau_{s-1}^* \} \stackrel {p}{\rightarrow} 0$, for $1 \leq s
\leq r^*$, and from (\ref{cvl}) with $\tau_{s-1}^*+1 \leq \ell<
\tau_s^*$ that $P_m \{ \tau_s = \ell, \tau_{s-1}=\tau_{s-1}^* \}
\stackrel {p}{\rightarrow} 0$. \end{pf}

\section{Discussion and concluding remarks}\label{sec5}

The central limit theorem for particle filters in this paper and in
\cite{Cho04,Del04,DM08,Kun05} considers the case of fixed $T$ as the
number $m$ of particles becomes infinite. In addition to the standard
error estimates~(\ref{420}), one can use the following modification
that substitutes $\psi_T$ in $\hat\sigma{}^2(\psi _T)$ by
``out-of-sample'' estimates of $\psi_T$. Instead of generating
additional particles to accomplish this, we apply a sample-splitting
technique to the $m$ particles, which is similar to $k$-fold
cross-validation. The standard error estimates in the example in
Section~\ref{sec4.1} use this technique with $k=2$. Divide the $m$
particles into $k$ groups of equal size $r=\lfloor m/k \rfloor$ except
for the last one that may have a larger sample size. For definiteness
consider bootstrap resampling at every stage; the case of residual
resampling or occasional resampling can be treated similarly. Denote
the particles at the $t$th generation, before resampling, by $\{
\widetilde{\mathbf{X}}{}_t^{ij}\dvtx  1 \leq i \leq r, 1 \leq j \leq k
\}$, and let $\mathbf{X}_t^{1j}, \ldots, \mathbf {X}_t^{rj}$ be sampled
with replacement from $\{ \widetilde{\mathbf{X}}{}_t^{1j}, \ldots,
\widetilde{\mathbf{X}}{}_t^{rj} \}$, with weights $W_t^{ij} =
w_t(\widetilde{\mathbf{X}}{}_t^{ij})/(m \bar w_t^j)$, where $\bar w_t^j
= r^{-1} \sum_{i=1}^r w_t(\widetilde{\mathbf{X}}{}_t^{ij})$. With these
``stratified'' resampling weights, we estimate $\psi_T$ by
%
\begin{equation}
\label{tpsiT} \hat\psi_T = k^{-1} \sum
_{j=1}^k \hat\psi_T^j\qquad\mbox{where } \hat\psi_T^j = \bigl(r \bar
w_T^j \bigr)^{-1} \sum
_{i=1}^r \psi \bigl(\widetilde{\mathbf{X}}{}_T^{ij}
\bigr) w_T \bigl(\widetilde{\mathbf{X}}{}_T^{ij}
\bigr).
\end{equation}
Similarly, letting $\hat\psi_T^{-j} = (k-1)^{-1} \sum_{\ell\dvtx
\ell\neq j} \hat\psi_T^{\ell}$, we estimate $\sigma^2$ by
%
\begin{equation}
\label{sSP} \hat\sigma{}^2_{\mathrm{SP}} = m^{-1} \sum
_{j=1}^k \sum_{i=1}^r
\biggl\{ \sum_{\ell\dvtx  A_{T-1}^{\ell j}
=i} \frac{w_T(\widetilde{\mathbf{X}}{}_T^{\ell j})}{\bar w_T^j} \bigl[\psi
\bigl(\widetilde{\mathbf{X}}{}_T^{\ell j} \bigr)-\hat
\psi_T^{-j} \bigr] \biggr\}^2.
\end{equation}

Chopin (Section~2.1 of \cite{Cho04}) summarizes a general framework of
traditional particle filter estimates of $\psi_T$, in which the
resampling (also called ``selection'') weights $w_t(\mathbf{X}_t)$ are
chosen to convert a weighted approximation of the posterior density of
$\mathbf{X}_t$ given $\mathbf{Y}_t$, with likelihood ratio weights
associated with importance sampling (called ``mutation'' in
\cite{Cho04}), to an unweighted approximation (that has equal weights
1), so that the usual average $\bar\psi_T:= m^{-1} \sum_{i=1}^m
\psi(\mathbf {X}_T^{i})$ converges a.s. to $\psi_T$ as $m
\rightarrow\infty$. He uses induction on $t \leq T$ to prove the
central limit theorem for $\sqrt{m}(\bar\psi_T-\psi_T)$: ``conditional
on past iterations, each step generates independent (but not
identically distributed) particles, which follow some (conditional)
central limit theorem'' (page~2392 of \cite{Cho04}), and he points out that
the particle filter of Gilks and Berzuini \cite{GB01} is a variant of
this framework.

Let $V_t$ be the variance of the limiting normal distribution of $\sqrt
{m_t}(\bar\psi_T-\psi_T)$ in the Gilks--Berzuini particle filter
(pages~132--134 of \cite{GB01}), that assumes
%
\begin{equation}
\label{m1t} m_1 \rightarrow\infty,\qquad\mbox{then } m_2
\rightarrow\infty, \ldots \mbox{ and then } m_T \rightarrow\infty,
\end{equation}
which means that $m_1 \rightarrow\infty$ and $m_t -m_{t-1}
\rightarrow\infty$ for $1 < t \leq T$. The estimates of $V_t$ proposed
in \cite{GB01} use the idea of ``ancestor'' instead of the ancestral
origin we use. For $s<t$, particle $\widetilde{\mathbf{X}}{}_s^{i}$ is
called an ancestor of $\widetilde{\mathbf{X}}{}_t^{k}$ if
$\widetilde{\mathbf{X}}{}_t^{k}$ descends from~$\widetilde{\mathbf{X}}{}_s^{i}$. Thus, the ancestral origin is the
special case corresponding to $s=1$ (the first generation of the
particles). Let
%
\begin{eqnarray}
\label{Ntkl} N_t^{k,\ell} & = & \sum
_{s=1}^t \sum_{i=1}^{m_s}
\mathbf{1}_{\{
\widetilde{\mathbf{X}}{}_s^{i}  \mathrm{\ is\ an\  ancestor\
of\  } \widetilde{\mathbf{X}}{}_t^{k} \mathrm{\ and\  of\  } \widetilde{\mathbf{X}}{}_t^{\ell} \}},
\\
\label{Vt} \widehat V_t & = & \frac{1}{m_t^2} \sum
_{k=1}^{m_t} \sum_{\ell
=1}^{m_t}
N_t^{k,\ell} \bigl\{ \psi \bigl(\widetilde{\mathbf{X}}{}_t^{k} \bigr)-\bar\psi_T \bigr\}
\bigl\{ \psi \bigl(\widetilde{\mathbf{X}}{}_t^{\ell} \bigr)-
\bar\psi_T \bigr\}.
\end{eqnarray}
Theorem 3 of Gilks and Berzuini (Appendix~A of \cite{GB01}) shows that
$\widehat V_t/V_t \stackrel{p}{\rightarrow} 1$ under (\ref{m1t}). The
basic idea is to use the law of large numbers for each generation
conditioned on the preceding one and an induction argument which relies
on the assumption (\ref{m1t}) that ``is not directly relevant to the
practical context.'' In contrast, our Theorem \ref{thm3a} or
\ref{thm4a} is based on a much more precise martingale approximation of
$m(\hat\psi_T-\psi_T)$ or $m(\hat\psi_{\mathrm{OR}}-\psi_T)$. While we
have focused on importance sampling in this paper, another approach to
choosing the proposal distribution is to use Markov chain Monte Carlo
iterations for the mutation step, as in \cite{GB01}, and important
advances in this approach have been developed recently in
\cite{ADH10,BDD10}.

\begin{appendix}\label{app}
\section*{\texorpdfstring{Appendix: Lemmas \lowercase{\protect\ref{lem1}, \protect\ref{lem1b},
\protect\ref{lem7}}, Corollary \lowercase{\protect\ref{thm2}} and
(\lowercase{\protect\ref{eq1}})--(\lowercase{\protect\ref{appendix4}})}
{Appendix: Lemmas 2, 3, 5, Corollary 1 and (3.30)--(3.32)}}\label{sec6}

\setcounter{equation}{0}

\begin{pf*}{Proof of Lemma \ref{lem1}}
We use induction on $t$ and show
that if (ii) holds for $t-1$, then (i) holds for $t$, and that if (i)
holds for $t$, then (ii) also holds for $t$. Since $G=G^+ - G^-$, we
can assume without loss of generality that $G$ is nonnegative. By the
law of large numbers, (i) holds for $t=1$, noting that $h_0^* \equiv1$.
Let $t>1$ and assume that (ii) holds for $t-1$. To show that (i) holds
for $t$, assume $\mu_t = E_q
[G(\mathbf{X}_t)/h_{t-1}^*(\mathbf{X}_{t-1})]$ to be finite and suppose
that, contrary to (\ref{lem11}), there exist $0 < \varepsilon< 1$ and
$m_1 < m_2 < \cdots$ such that
%
\begin{equation}
\label{aA2} \qquad P_m \Biggl\{ \Biggl| m^{-1} \sum
_{i=1}^m G \bigl(\widetilde{\mathbf{X}}{}_{t,m}^i
\bigr) - \mu_t \Biggr| > 2 \varepsilon \Biggr\} > \varepsilon(2+
\mu_t)\qquad\mbox{for all } m \in\mathcal{M},
\end{equation}
in which $\mathcal{M}= \{ m_1, m_2, \ldots\}$ and we write
$\widetilde{\mathbf{X}}{}_{t,m}^i (=\widetilde{\mathbf{X}}{}_t^i)$ to
highlight the dependence on $m$. Let $\delta= \varepsilon^3$. Since
$x=(x \wedge y)+(x-y)^+$, we can write
$G(\widetilde{\mathbf{X}}{}_{t,m}^i) = U_{t,m}^i + V_{t,m}^i +
S_{t,m}^i$, where
%
\begin{eqnarray}
\qquad U_{t,m}^i & = & G \bigl(\widetilde{\mathbf{X}}{}_{t,m}^i \bigr) \wedge(\delta m)-E_m
\bigl[G \bigl(\widetilde{\mathbf{X}}{}_{t,m}^i \bigr) \wedge(
\delta m)|\mathcal{F}_{2t-2} \bigr],
\nonumber\\[-8pt]\label{UVY}  \\[-8pt]
V_{t,m}^i & = & E_m \bigl[G \bigl(
\widetilde{\mathbf{X}}{}_{t,m}^i \bigr) \wedge(\delta m)|
\mathcal{F}_{2t-2} \bigr], \qquad S_{t,m}^i =
\bigl[G \bigl(\widetilde{\mathbf{X}}{}_{t,m}^i \bigr) -\delta
m \bigr]^+.\nonumber
\end{eqnarray}
Since $E_m(U_{t,m}^i|\mathcal{F}_{2t-2})=0$ and $\operatorname
{Cov}_m(U_{t,m}^i,U_{t,m}^{\ell}| \mathcal{F}_{2t-2})=0$ for $i
\neq\ell$, it follows from Chebyshev's inequality, $(U_{t,m}^i)^2 \leq
\delta m G(\widetilde{\mathbf{X}}{}_{t,m}^i)$ and $\delta=
\varepsilon^3$ that
%
\begin{eqnarray}
\label{var} && P_m \Biggl\{ \Biggl| m^{-1} \sum
_{i=1}^m U_{t,m}^i \Biggr| >
\varepsilon \bigg| \mathcal{F}_{2t-2} \Biggr\}\nonumber
\\
&&\qquad \leq (\varepsilon m)^{-2} \operatorname{Var}_m \Biggl( \sum_{i=1}^m U_{t,m}^i \bigg| \mathcal{F}_{2t-2} \Biggr)
\nonumber\\[-8pt]\\[-8pt]
&&\qquad =  (\varepsilon m)^{-2} \sum_{i=1}^m E_m \bigl(\bigl|U_{t,m}^i\bigr|^2|
\mathcal{F}_{2t-2} \bigr)\nonumber
\\
&&\qquad \leq\varepsilon m^{-1} \sum
_{i=1}^m E_m \bigl[G \bigl(\widetilde{\mathbf{X}}{}_{t,m}^i \bigr)|\mathcal{F}_{2t-2}
\bigr] \stackrel{p} {\rightarrow} \varepsilon\mu_t,
\nonumber
\end{eqnarray}
as $m \rightarrow\infty$, by (\ref{lem13}) applied to $G^*(\mathbf
{X}_{t-1}) = E_q[ G(\mathbf{X}_t)|\mathbf{X}_{t-1}]$, noting that
$E_q[G^*(\mathbf {X}_{t-1})/h_{t-1}^*(\mathbf{X}_{t-1})] = \mu_t$.
Application of (\ref{lem13}) to $G^*(\mathbf{X}_{t-1})$ also yields
%
\begin{equation}
\label{mkVk} m^{-1} \sum_{i=1}^m
V_{t,m}^i \stackrel{p} {\rightarrow} \mu_t.
\end{equation}
From (\ref{var}), it follows that
%
\begin{equation}
\label{aA4} P_m \Biggl\{ \Biggl| m^{-1} \sum
_{i=1}^m U_{t,m}^i \Biggr| >
\varepsilon \Biggr\} \leq\varepsilon(1+\mu_t)\qquad \mbox{for all large } m.
\end{equation}
Since $\mu_t < \infty$, we can choose $n_k$ such that
\[
E_q \bigl[G(\mathbf{X}_t) \mathbf{1}_{\{ G(\mathbf{X}_t) > n_k \}
}/h_{t-1}^*(
\mathbf{X}_{t-1}) \bigr] < k^{-1}.
\]
Hence, by applying (\ref{lem13}) to $G_k(\mathbf{X}_{t-1}) =
E_q[G(\mathbf{X}_t) \mathbf{1}_{\{ G(\mathbf{X}_t)>n_k
\}}|\mathbf{X}_{t-1}]$, we can choose a further subsequence $m_k$ that
satisfies (\ref{aA2}), $m_k \geq n_k/\delta$ and
\[
P_m \Biggl\{ m^{-1} \sum_{i=1}^m
E_m \bigl[G \bigl(\widetilde{\mathbf{X}}{}_{t,m}^i
\bigr) \mathbf{1}_{\{ G(\widetilde{\mathbf{X}}{}_{t,m}^i) >n_k \}}|\mathcal {F}_{2t-2} \bigr]
\geq2k^{-1} \Biggr\} \leq k^{-1}\quad\mbox{for }
m=m_k. 
\]
Then for $m=m_k$, with probability at least $1-k^{-1}$,
\begin{eqnarray*}
\sum_{i=1}^m P_m \bigl\{
S_{t,m}^i \neq 0 | \mathcal{F}_{2t-2} \bigr\} &=&
\sum_{i=1}^m P_m \bigl\{ G
\bigl(\widetilde{\mathbf{X}}{}_{t,m}^i \bigr) > \delta m|
\mathcal{F}_{2t-2} \bigr\}
\\
&\leq& (\delta m)^{-1} \sum
_{i=1}^m E_m \bigl[ G \bigl(
\widetilde{\mathbf{X}}{}_{t,m}^i \bigr) \mathbf{1}_{\{ G(\widetilde{\mathbf{X}}{}_{t,m}^i) > \delta m \}}
| \mathcal{F}_{2t-2} \bigr]
\\
&\leq& 2 \delta^{-1}
k^{-1}.
\end{eqnarray*}
Therefore, $P_m \{ S_{t,m}^i \neq0$ for some $1 \leq i \leq m \}
\rightarrow0$ as $m=m_k \rightarrow\infty$. Combining this with
(\ref{mkVk}) and (\ref{aA4}), we have a
contradiction to (\ref{aA2}). Hence, we conclude (\ref{lem11}) holds
for $t$. The proof that (ii) holds for $t$ whenever (i) holds for $t$
is similar. In particular, since (i) holds for $t=1$, so does (ii).
\end{pf*}

\begin{pf*}{Proof of Lemma \ref{lem1b}}
We again use
induction and show that if (ii) holds for $t-1$, then (i) holds for
$t$, and that if (i) holds for $t$, then so does (ii). First (i) holds
for $t=1$ because $A_0^i =i$, $h_0^* \equiv1$ and $E_q G(X_1) < \infty$
implies that for any $\varepsilon> 0$,
\[
P_m \Bigl\{ m^{-1} \max_{1 \leq i \leq m} G \bigl(
\widetilde{X}{}_1^i \bigr) \geq\varepsilon \Bigr\} 
\leq m P_q \bigl\{ G(X_1) \geq
\varepsilon m \bigr\} \rightarrow0\quad\mbox{as } m \rightarrow \infty.
\]
Let $t > 1$ and assume that (ii) holds for $t-1$. To show that (i)
holds for $t$, suppose $\mu_t (=E_q[G(\mathbf{X}_t)/
h_{t-1}^*(\mathbf{X}_{t-1})]) < \infty$ and that, contrary to
(\ref{lem21}), there exist $0 < \varepsilon< 1$ and $\mathcal{M}= \{
m_1, m_2, \ldots\} $ with $m_1 < m_2 < \cdots$ such that
%
\begin{equation}
\label{i2a}  P_m \biggl\{ m^{-1} \max
_{1 \leq j \leq m} \sum_{i\dvtx A_{t-1}^i =j} G \bigl(
\widetilde{\mathbf{X}}{}_{t,m}^i \bigr) > 2 \varepsilon \biggr\}
> \varepsilon(2+\mu_t)\quad\mbox{for all } m \in\mathcal{M}.\hspace*{-40pt}
\end{equation}
As in (\ref{UVY}), let $\delta=\varepsilon^3$ and
$G(\widetilde{\mathbf{X}}{}_{t,m}^i) = U_{t,m}^i+V_{t,m}^i +S_{t,m}^i$.
For $1 \leq j \leq m$,
\[
P_m \biggl\{ \biggl| m^{-1} \sum_{i\dvtx A_{t-1}^i=j}
U_{t,m}^i \biggr| > \varepsilon \Big| \mathcal{F}_{2t-2} \biggr
\} \leq\varepsilon m^{-1} \sum_{i\dvtx A_{t-1}^i=j}
E_m \bigl[G \bigl(\widetilde{\mathbf{X}}{}_{t,m}^i
\bigr)|\mathcal{F}_{2t-2} \bigr],
\]
by an argument similar to (\ref{var}). Summing over $j$ and then taking
expectations, it then follows from Lemma \ref{lem1} that
%
\begin{equation}
\label{mU} \qquad\sum_{j=1}^m
P_m \biggl\{ \biggl| m^{-1} \sum_{i\dvtx A_{t-1}^i=j}
U_{t,k}^i \biggr| > \varepsilon \biggr\} \leq\varepsilon(1+
\mu_t)\quad\mbox{for all large } m.
\end{equation}
By (\ref{lem22}) applied to the function $G^*(\mathbf{X}_{t-1}) =
E_q[G(\mathbf{X}_t)|\mathbf{X}_{t-1}]$,
%
\begin{equation}
\label{mV} m^{-1} \max_{1 \leq j \leq m} \sum
_{i\dvtx  A_{t-1}^i =j} V_{t,m}^i \stackrel{p} {
\rightarrow} 0\qquad\mbox{as } m \rightarrow\infty.
\end{equation}
As in the proof of Lemma \ref{lem1}, select a further subsequence
$m_k'$ such that (\ref{i2a}) holds and $P_m \{ \sum_{i=1}^m S_{t,m}^i =
0 \} \rightarrow1$ as $m=m_k' \rightarrow\infty$. Combining this with
(\ref{mU}) and (\ref{mV}) yields a contradiction to (\ref{i2a}). Hence,
(\ref{lem21}) holds for~$t$.

Next let $t \geq1$ and assume that (i) holds for $t$. To show that (ii)
holds for $t$, assume $\tilde\mu_t =
E_q[G(\mathbf{X}_t)/h_t^*(\mathbf{X}_t)]$ to be\vspace*{-2pt} finite and note that
$\sum_{i\dvtx A_t^i=j} G(\mathbf{X}_{t,k}^i) = \sum_{i\dvtx
A_{t-1}^i=j} \# _{t,k}^i G(\widetilde{\mathbf{X}}{}_{t,k}^i)$. Suppose
that, contrary to (\ref{lem22}), there exist $0 < \varepsilon< 1$ and
$\mathcal{M}= \{ m_1, m_2, \ldots\}$ with $m_1 < m_2 < \cdots$ such
that
%
\begin{eqnarray}
P_m \biggl\{ m^{-1} \max
_{1 \leq j \leq m} \sum_{i\dvtx  A_{t-1}^i =j}
\#_{t,k}^i G \bigl(\widetilde{\mathbf{X}}{}_{t,k}^i
\bigr) > 2 \varepsilon \biggr\} > \varepsilon(2+\tilde\mu_t)
\nonumber\\[-18pt]\label{wtu} \\[-6pt]
\eqntext{\mbox{for all } m \in\mathcal{M}.}
\end{eqnarray}
Let $\delta=\varepsilon^3$ and write $\#_{t,m}^i
G(\widetilde{\mathbf{X}}{}_{t,m}^i) = \widetilde U_{t,m}^i + \widetilde
V_{t,m}^i + \#_{t,m}^i S_{t,m}^i$, where
\begin{eqnarray*}
\widetilde U_{t,m}^i &=& \bigl(\#_{t,m}^i-m
W_{t,m}^i \bigr) \bigl[G \bigl(\widetilde{\mathbf{X}}{}_{t,m}^i \bigr) \wedge(\delta m) \bigr],
\\
\widetilde V_{t,m}^i &=& m W_{t,m}^i
\bigl[G \bigl(\widetilde{\mathbf{X}}{}_{t,m}^i \bigr) \wedge(
\delta m) \bigr]
\end{eqnarray*}
and $S_{t,m}^i$ is defined in (\ref{UVY}). Since $E_m(\#_{t,m}^i|
\mathcal{F}_{2t-1}) = m W_{t,m}^i$,
$\operatorname{Var}_m(\#_{t,m}^i|\break
\mathcal{F}_{2t-1}) \leq m W_{t,m}^i =
w_t(\widetilde{\mathbf{X}}{}_{t,m}^i)/\bar w_{t,m}$ and since
$\#_{t,m}^1, \ldots, \#_{t,m}^m$ are pairwise negatively correlated
conditioned on $\mathcal {F}_{2t-1}$, we obtain by Chebyshev's
inequality that for $1 \leq j \leq m$, $P_m \{ | m^{-1} \sum_{i\dvtx
A_{t-1}^i=j} \widetilde U_{t,m}^i | > \varepsilon| \mathcal{F}_{2t-1}
\}$ is bounded above by
\begin{eqnarray*}
&& \frac{1}{(\varepsilon m)^{2}} \sum_{i\dvtx A_{t-1}^i=j} \operatorname{Var}_m
\bigl( \#_{t,m}^i | \mathcal{F}_{2t-1} \bigr)
\bigl[G \bigl(\widetilde{\mathbf{X}}{}_{t,m}^i \bigr) \wedge(
\delta m) \bigr]^2
\\
&&\qquad \leq \frac{\varepsilon}{m \bar w_{t,m}} \sum
_{i\dvtx A_{t-1}^i=j} G^* \bigl(\widetilde{\mathbf{X}}{}_{t,m}^i
\bigr),
\end{eqnarray*}
where $G^*(\cdot) = w_t(\cdot) G(\cdot)$. Hence, applying Lemma
\ref{lem1} to $G^*$ and noting that, by~(\ref{etah}),
$E_q[G^*(\mathbf{X}_t)/h_{t-1}^*(\mathbf{X}_{t-1})] =\zeta_t^{-1}
\tilde\mu_t$ and $\bar w_{t,m} = m^{-1} \sum_{i=1}^m
w_m(\widetilde{\mathbf{X}}{}_{t,m}^i) \stackrel{p}{\rightarrow}
\zeta_t^{-1}$ by Lemma \ref{lem1}, we obtain
%
\begin{equation}
\label{mtU} \qquad\sum_{j=1}^m
P_m \biggl\{ \biggl| m^{-1} \sum_{i\dvtx A_{t-1}^i=j}
\widetilde U_{t,m}^i \biggr| > \varepsilon \biggr\} \leq\varepsilon(1+
\tilde\mu_t)\qquad\mbox{for all large } m.
\end{equation}
By (\ref{lem21}) applied to $G^*$ and Lemma \ref{lem1},
%
\begin{equation}
\label{mtV} m^{-1} \max_{1 \leq j \leq m} \sum
_{i\dvtx A_{t-1}^i=j} \widetilde V_{t,m}^i \leq \bar
w_{t,m}^{-1} \biggl[ m^{-1} \max
_{1 \leq j \leq m} \sum_{i\dvtx A_{t-1}^i=j} G^* \bigl(
\widetilde{\mathbf{X}}{}_{t,m}^i \bigr) \biggr] \stackrel{p} {
\rightarrow} 0.\hspace*{-40pt}
\end{equation}
Finally, by Lemma \ref{lem1} applied to $G_m(\mathbf{X}_t) =
w_t(\mathbf{X}_t) G(\mathbf{X}_t) \mathbf{1}_{\{ G(\mathbf{X}_t) >
\delta m \}}$,
%
\begin{eqnarray}
\label{A14}
&& \sum_{i=1}^m
P_m \bigl\{ \#_{t,m}^i S_{t,m}^i
\neq0 | \mathcal {F}_{2t-1} \bigr\}\nonumber
\\
&&\qquad = \sum_{i=1}^m P_m \bigl\{ \#_{t,m}^i > 0 |
\mathcal{F}_{2t-1} \bigr\} \mathbf {1}_{\{ S_{t,m}^i \neq0 \}}
\nonumber\\[-8pt]\\[-8pt]
&&\qquad \leq\frac{1}{\delta m} \sum_{i=1}^m
E_m \bigl(\#_{t,m}^i|\mathcal{F}_{2t-1}
\bigr) \frac{G_m(\widetilde{\mathbf{X}}{}_{t,m}^i)} 
{w_t (\widetilde{\mathbf{X}}{}_{t,m}^i )}\nonumber
\\
&&\qquad = \frac{1}{\delta\bar w_{t,m}}
\Biggl[ \frac{1}{m} \sum_{i=1}^m
G_m \bigl(\widetilde{\mathbf{X}}{}_{t,m}^i \bigr)
\Biggr] \stackrel{p} {\rightarrow} 0 \nonumber
\end{eqnarray}
as $m=m_k \rightarrow\infty$. Since (\ref{mtU})--(\ref{A14})
contradicts (\ref{wtu}), (\ref{lem22}) follows.
\end{pf*}

\begin{pf*}{Proof of Lemma \ref{lem7}}
The proof of Lemma \ref{lem7} is
similar to that of Lemma~\ref{lem1b}, and in fact is simpler because
the finiteness assumption on the second moments of $G$ dispenses with
truncation arguments of the type in (\ref{UVY}).
First (i) holds for $t=1$ because $A_0^i=i$, $h_0^*=1$ and
%
\[
P_m \Biggl\{ m^{-1} \sum_{i=1}^m
G^2 \bigl(\widetilde{X}{}_1^i \bigr) \geq K
\Biggr\} \leq K^{-1} E_q G^2(X_1)
\rightarrow0\qquad\mbox{as } K \rightarrow\infty.
\]
As in the proof of Lemma \ref{lem1b}, let $t>1$ and assume that (ii)
holds for $t-1$. To show that (i) holds for $t$, suppose
$\mu_t:=E_q[G^2(\mathbf{X}_t) \Gamma_{t-1}(\mathbf{X}_{t-1})] <
\infty$. By expressing $G=G^+-G^-$, we may assume without loss of
generality $G$ is nonnegative. Write $G(\widetilde{\mathbf{X}}{}_t^i) =
U_t^i + V_t^i$, where $V_t^i =
E_m[G(\widetilde{\mathbf{X}}{}_t^i)|\mathbf{X}_{t-1}^i]$. By the
induction hypothesis, $m^{-1} \sum_{j=1}^m ( \sum_{i\dvtx A_{t-1}^i=j}
V_t^i)^2 = O_p(1)$. To show\vspace*{-1pt} that $m^{-1} \sum_{j=1}^m (\sum_{i\dvtx
A_{t-1}^i=j} U_t^i)^2 = O_p(1)$, note that $E_m(U_t^i
U_t^\ell|\mathcal{F}_{2t-2})=0$ for $i \neq\ell$. Hence, it follows by
an argument similar to (\ref{var}) that
\begin{eqnarray*}
&& P_m \Biggl\{ m^{-1} \sum_{j=1}^m
\biggl( \sum_{i\dvtx A_{t-1}^i=j} U_t^i
\biggr)^2 \geq K \Big| \mathcal{F}_{2t-2} \Biggr\}
\\
&&\qquad \leq(Km)^{-1} \sum_{i=1}^m
E_m \bigl[ \bigl(U_t^i \bigr)^2|
\mathcal{F}_{2t-2} \bigr]
\\
&&\qquad \stackrel{p} {\rightarrow}
K^{-1} E_q \bigl\{ \bigl|G(\mathbf{X}_t)-E_q
\bigl[G(\mathbf{X}_t)|\mathbf {X}_{t-1}
\bigr]\bigr|^2/h_{t-1}^*( \mathbf{X}_{t-1}) \bigr\}.
\end{eqnarray*}

Next assume that (i) holds for $t$ and show that (ii) holds for $t$ if
$\tilde\mu_t:=E_q[G^2(\mathbf{X}_t) \Gamma_t(\mathbf{X}_t)]$ is finite.
Note that $\sum_{i\dvtx A_t^i=j} G(\mathbf{X}_t^i) = \sum_{i\dvtx
A_{t-1}^i=j} \#_t^i G(\widetilde{\mathbf{X}}{}_t^i)$. Write $\#_t^i
G(\widetilde{\mathbf{X}}{}_t^i) = \widetilde U_t^i + \widetilde V_t^i$,
where $\widetilde V_t^i = mW_t^i G(\widetilde{\mathbf{X}}{}_t^i) =
w_t(\widetilde{\mathbf{X}}{}_t^i) G(\widetilde{\mathbf{X}}{}_t^i)/\bar
w_t$. Since $\bar w_t \stackrel{p}{\rightarrow} \zeta_t^{-1}$ and since
$E_q \{ [w_t(\mathbf{X}_t) G(\mathbf{X}_t)]^2
\Gamma_{t-1}(\mathbf{X}_{t-1}) \} \leq\tilde\mu_t < \infty$, it follows
from the induction hypothesis that $m^{-1} \sum_{j=1}^m ( \sum_{i\dvtx
A_{t-1}^i=j} \widetilde V_t^i)^2 = O_p(1)$.\vspace*{-1pt} Moreover, since
$\operatorname{Var}_m(\#_t^i|\mathcal{F}_{2t-1}) \leq mW_t^i$ and
$\operatorname{Cov}_m(\#_t^i,\#_t^\ell|\mathcal{F}_{2t-1}) \leq0$ for
$i \neq\ell$, it follows from an argument similar to (\ref{var}) that
\begin{eqnarray*}
P_m \Biggl\{ m^{-1} \sum_{j=1}^m
\biggl( \sum_{i\dvtx A_{t-1}^i=j} \widetilde U_t^i
\biggr)^2 \geq K \Big| \mathcal{F}_{2t-1} \Biggr\}
& \leq& (Km)^{-1} \sum_{i=1}^m
mW_t^i G^2 \bigl(\widetilde{\mathbf{X}}{}_t^i \bigr)
\\
&\stackrel{p} {\rightarrow}& K^{-1} E_q \bigl[G^2(\mathbf
{X}_t)/h_t^*(\mathbf{X}_t) \bigr],
\end{eqnarray*}
noting that $E_q[G^2(\mathbf{X}_t)/h_t^*(\mathbf{X}_t)] \leq\eta
_t^{-1} \tilde\mu_t < \infty$. Hence,
\[
m^{-1} \sum_{j=1}^m \biggl(
\sum_{i\dvtx A_{t-1}^i=j} \widetilde U_t^i\biggr)^2 = O_p(1).
\]
This shows that (ii) holds for $t$, completing the induction proof.
\end{pf*}

\begin{pf*}{Proof of Corollary \ref{thm2}}
Recall that $\{(Z_k^{1},\ldots,Z_k^{m}), \mathcal{F}_k, 1 \leq k
\leq\break 2T-1 \}$ is a martingale difference sequence and that $Z_k^{1},
\ldots, Z_k^{m}$ are conditionally independent given
$\mathcal{F}_{k-1}$, where $\mathcal{F}_0$ is the trival
$\sigma$-algebra; moreover,
%
\begin{equation}
\label{sigb} \sqrt{m} (\tilde\psi_T - \psi_T) = \sum
_{k=1}^{2T-1} \Biggl( \sum
_{i=1}^m Z_k^{i}/\sqrt{m}
\Biggr).
\end{equation}
Fix $k$. Conditioned on $\mathcal{F}_{k-1}$,
$Y_{mi}^{k}:=Z_k^{i}/\sqrt{m}$ is a triangular array of row-wise
independent (i.i.d. when $k=2t$) random variables. Noting that
$E_m[\tilde f_t(\widetilde{\mathbf{X}}{}_t^{i})|\break
\mathbf{X}_{t-1}^{i}]=\tilde f_{t-1}(\mathbf{X}_{t-1}^{i})$, we obtain
from (\ref{list}) that
\begin{eqnarray*}
E_m \bigl[ \bigl(Z_{2t-1}^{i}
\bigr)^2|\mathcal{F}_{2t-2} \bigr] & = & \bigl\{
E_m \bigl[\tilde f_t^2 \bigl(\widetilde{\mathbf{X}}{}_t^{i} \bigr)|\mathbf {X}_{t-1}^{i}
\bigr]-\tilde f^2_{t-1} \bigl(\mathbf{X}_{t-1}^{i}
\bigr) \bigr\} \bigl(H_{t-1}^i \bigr)^2,
\\
E_m \bigl[ \bigl(Z_{2t}^{i}
\bigr)^2|\mathcal{F}_{2t-1} \bigr] & = & \sum
_{j=1}^m W_t^{j} \tilde
f_t^2 \bigl(\widetilde{\mathbf{X}}{}_t^{j}
\bigr) \bigl(\widetilde{H}{}_t^j \bigr)^2 -
\Biggl[\sum_{j=1}^m W_t^{j}
\tilde f_t \bigl(\widetilde{\mathbf{X}}{}_t^{j}
\bigr) \widetilde{H}{}_t^j \Biggr]^2.
\end{eqnarray*}
Let $\varepsilon> 0$. Since
$W_t^{j}=w_t(\widetilde{\mathbf{X}}{}_t^{j})/(m \bar w_t)$ and $\bar
w_t \stackrel{p}{\rightarrow} \zeta_t^{-1}$,
\[
\sum_{i=1}^m E_m \bigl[
\bigl(Y_{mi}^{k} \bigr)^2|\mathcal{F}_{k-1}
\bigr] \stackrel {p} {\rightarrow} \tilde\sigma_k^2,
\qquad\sum_{i=1}^m
E_m \bigl[ \bigl(Y_{mi}^{k}
\bigr)^2 \mathbf{1}_{\{ |Y_{mi}^{k}| >
\varepsilon\}} |\mathcal{F}_{k-1} \bigr]
\stackrel{p} { \rightarrow} 0,
\]
by (\ref{etah}) and Lemma~\ref{lem1}. Hence, by Lindeberg's central
limit theorem for triangular arrays of independent random variables,
the conditional distribution of $\sum_{i=1}^m Z_k^{i}/\sqrt{m}$ given
$\mathcal {F}_{k-1}$ converges to $N(0,\tilde\sigma_k^2)$ as $m
\rightarrow\infty$. This implies that for any real $u$,
%
\begin{equation}
\label{A20} E_m \bigl(e^{\mathbf{i} u \sum_{j=1}^m Z_k^j/\sqrt{m}}|\mathcal{F}_{k-1}
\bigr) \stackrel{p} {\rightarrow} e^{-u^2 \tilde\sigma_k^2/2}, \qquad1 \leq k \leq2T-1,
\end{equation}
where $\mathbf{i}$ denotes the\vspace*{-2pt} imaginary number $\sqrt{-1}$. Let $S_t =
\sum_{k=1}^t (\sum_{j=1}^m Z_k^j/ \sqrt{m})$. Since $\{ \sum_{j=1}^m
Z_k^j/\sqrt{m}$, $\mathcal{F}_k,  1 \leq k \leq2T-1 \}$ is a martingale
difference sequence, it follows from (\ref{sigb}), (\ref{A20}) and
mathematical induction that
\begin{eqnarray*}
&& E_m \bigl(e^{\mathbf{i} u \sqrt{m}(\tilde\psi_T
-\psi_T)} \bigr)
\\
&&\qquad = E_m \bigl[e^{\mathbf{i} uS_{2T-2}} E_m \bigl(e^{\mathbf{i} u \sum
_{j=1}^m Z_{2T-1}^j/\sqrt{m}} |
\mathcal{F}_{2T-2} \bigr) \bigr]
\\
&&\qquad =  e^{-u^2 \tilde\sigma_{2T-1}^2/2}
E_m \bigl(e^{\mathbf{i}
uS_{2T-2}} \bigr)+o(1) = \cdots= \exp
\bigl(-u^2 \sigma^2_{\mathrm{C}}/2 \bigr)+o(1).
\end{eqnarray*}
Hence, (\ref{clt}) holds.
\end{pf*}

\begin{pf*}{Proof of (\ref{eq1})--(\ref{appendix4})}
$\!\!\!$For distinct $i$, $j$,
$k$, $\ell$, define $c_{ij}=(mW_t^i)(mW_t^j)$, $c_{ijk \ell} =
(mW_t^i)(mW_t^j)(mW_t^k)(mW_t^{\ell})$, $c_{iijk} = (mW_t^i)^2
(mW_t^j)(mW_t^k)$ and define $c_{ijk}$, $c_{iij}$, $c_{ijj}$ similarly.
Since $\#_t^i = \sum_{h=1}^m \mathbf{1}_{\{ B_t^h =i \}}$, where
$B_t^h$ are i.i.d. conditioned on $\mathcal{F}_{2t-1}$ such that $P \{
B_t^h=i | \mathcal{F}_{2t-1} \} = W_t^i$, combinatorial arguments show
that
%
\begin{eqnarray}
\qquad\quad && E_m \bigl[ \bigl(\#_t^i-mW_t^i
\bigr) \bigl(\#_t^j-mW_t^j \bigr)
\bigl(\#_t^k-mW_t^k \bigr) \bigl(
\#_t^{\ell}-mW_t^{\ell} \bigr) |
\mathcal{F}_{2t-1} \bigr]\nonumber
\\
&&\qquad = m^{-4} \bigl[m(m-1) (m-2) (m-3)
\nonumber\\[-8pt]\label{238a} \\[-8pt]
&&\hspace*{54pt}{}-4m^2(m-1)
(m-2)+6m^3(m-1)-4m^4+m^4 \bigr] c_{ijk \ell}\nonumber
\\
&&\qquad = \bigl(3m^{-2}-6m^{-3} \bigr) c_{ijk \ell}, \nonumber
\\
&& E_m \bigl[ \bigl(\#_t^i-mW_t^i
\bigr)^2 \bigl(\#_t^j-mW_t^j
\bigr) \bigl(\#_t^k-mW_t^k
\bigr)| \mathcal{F}_{2t-1} \bigr]
\nonumber\\[-8pt]\label{238b} \\[-8pt]
&&\qquad= \bigl(3m^{-2}-6m^{-3} \bigr) c_{iijk} +
\bigl(-m^{-1}+2m^{-2} \bigr) c_{ijk},\nonumber
\\
&& E_m \bigl[ \bigl(\#_t^i-mW_t^i
\bigr)^2 \bigl(\#_t^j-mW_t^j
\bigr)^2 | \mathcal {F}_{2t-1} \bigr]\nonumber
\\
\label{238c} &&\qquad = \bigl(3m^{-2}-6m^{-3} \bigr)c_{iijj} +
\bigl(-m^{-1}+2m^{-2} \bigr) (c_{iij}+c_{ijj})
\\
&&\quad\qquad{} + \bigl(1-m^{-1} \bigr) c_{ij}, \nonumber
\\
&& E_m \bigl[ \bigl(\#_t^i-mW_t^i\bigr)^4|\mathcal{F}_{2t-1} \bigr]\nonumber
\\
&&\qquad  = m \bigl[W_t^i\bigl(1-W_t^i \bigr)^4 +\bigl(1-W_t^i \bigr) \bigl(W_t^i\bigr)^4 \bigr]
\nonumber\\[-8pt]\label{4thmoment} \\[-8pt]
&&\quad\qquad{} + 3m(m-1) \bigl[ \bigl(W_t^i \bigr)^2 \bigl(1-W_t^i \bigr)^4\nonumber
\\
&&\hspace*{97pt}{} +2 \bigl(W_t^i \bigr)^3 \bigl(1-W_t^i \bigr)^3 + \bigl(W_t^i \bigr)^4 \bigl(1-W_t^i \bigr)^2 \bigr],\nonumber
\\
\label{cij} && E_m \bigl\{ \bigl(\#_t^i-mW_t^i \bigr) \bigl(\#_t^j-mW_t^j \bigr)| \mathcal{F}_{2t-1} \bigr\} = -m^{-1} c_{ij}.
\end{eqnarray}

To prove (\ref{eq1}), let $T_i = mW_t^i [\tilde
f_t(\widetilde{\mathbf{X}}{}_t^i) \widetilde{H}{}_t^i-\tilde f_0]^2$,
$\widebar T = m^{-1} \sum_{i=1}^m T_i$ and $T_t^* = m^{-1} \max_{1 \leq
j \leq m} \sum_{i\dvtx A_{t-1}^i=j} T_i$. We can use (\ref{238c}) and
(\ref{4thmoment}) to show that $m^{-2} E_m [ ( \sum_{i=1}^m S_i^2 )^2 |
\mathcal{F}_{2t-1} ]$ is equal to
\[
m^{-2} \sum_{i=1}^m
E_m \bigl(S_i^4|\mathcal{F}_{2t-1}
\bigr) +m^{-2} \sum_{i
\neq j} E_m
\bigl(S_i^2 S_j^2 |
\mathcal{F}_{2t-1} \bigr) = \widebar T{\,}^2 +o_P(1)
\]
as $m \rightarrow\infty$. By Lemma \ref{lem1}(i), $\widebar T{\,}^2
\stackrel{p}{\rightarrow} \tilde\sigma{}_{2t}^4$. From this\vspace*{1pt} and
(\ref{condV}), (\ref{eq1}) follows.

To prove (\ref{eq2}), let $D_i=mW_t^i [\tilde
f_t(\widetilde{\mathbf{X}}{}_t^i) \widetilde{H}{}_t^i-\tilde f_0]$,
$\widebar D = m^{-1} \sum_{i=1}^m |D_i|$ and $D_t^* = m^{-1} \max_{1
\leq j \leq m} \sum_{i\dvtx A_{t-1}^i=j} |D_i|$. By
(\ref{238a})--(\ref{238c}), the left-hand side of (\ref{eq2}) is equal
to
\begin{eqnarray*}
& & \frac{1}{m^2} \biggl( \frac{3}{m^2}-\frac{6}{m^3} \biggr)
\Biggl\{ \sum_{j=1}^m \sum
_{(i,\ell) \in C_{t-1}^j} D_i D_\ell \Biggr\}^2
+ \frac{1}{m^2} \biggl( 1-\frac{1}{m} \biggr) \sum
_{j=1}^m \sum_{(i,\ell) \in C_{t-1}^j}
T_i T_\ell
\\
& & \quad{}- \frac{4}{m^2} \biggl(
\frac{1}{m}-\frac{2}{m^2} \biggr) \sum_{j=1}^m
\biggl\{ \sum_{i\dvtx A_{t-1}^i =j} T_i \biggl( \sum
_{\ell\neq i\dvtx  A_{t-1}^{\ell}=j} D_\ell \biggr)^2 \biggr
\}
\\
&&\qquad \leq 3 \bigl(\widebar D D_t^* \bigr)^2 +
\widebar T T_t^*\qquad\mbox{for } m \geq2.
\end{eqnarray*}
By Lemmas \ref{lem1} and \ref{lem1b}, the upper bound in the above
inequality converges to 0 in probability as $m \rightarrow\infty$,
proving (\ref{eq2}).

By (\ref{cij}), $E_m(S_i S_j|\mathcal{F}_{2t-1}) = -m^{-1} D_i D_j$
and, similarly, $E_m(S_i^2|\mathcal{F}_{2t-1}) = T_i - m^{-1} D_i^2$.
Therefore, by (\ref{210}), the left-hand side of (\ref{appendix4}) is
equal to
\[
m^{-1} \sum_{j=1}^m \bigl(
\varepsilon_1^j+\cdots+\varepsilon_{2t-1}^j
\bigr)^2 \biggl\{ m^{-1} \sum_{i\dvtx A_{t-1}^i =j}
T_i - \biggl( m^{-1} \sum_{i\dvtx A_{t-1}^i
=j}
D_i \biggr)^2 \biggr\} \stackrel{p} {\rightarrow} 0,
\]
since $T_t^*+(D_t^*)^2 \stackrel{p}{\rightarrow} 0$ by Lemma
\ref{lem1b} and since (\ref{226}) holds for $k=2t-1$ by the induction
hypothesis.
\end{pf*}
\end{appendix}




\printaddresses

\end{document}